\documentclass[12 pt]{article}
\usepackage{algorithmic}
\usepackage{algorithm}
\usepackage{amsmath}
\usepackage{amssymb}
\usepackage{todonotes}
\newtheorem{theorem}{Theorem}
\newtheorem{remark}[theorem]{Remark}

\newtheorem{lemma}[theorem]{Lemma}

\newtheorem{proposition}[theorem]{Proposition}

\newtheorem{definition}[theorem]{Definition}

\newtheorem{corollary}[theorem]{Corollary}

\newtheorem{example}[theorem]{Example}

\def\urltilda{\kern -.15em\lower .7ex\hbox{\~{}}\kern .04em}

\date{} %

\newcommand{\ra}{\rightarrow}

\newcommand{\w}{{\rm w_H}}
\newcommand{\tr}{{\rm tr}}
\newcommand{\F}{{\mathbb{F}}}
\newcommand{\stab}{full degree-stability under restrictions to affine spaces}

\newcommand{\SqBinom}[2]{\genfrac{[}{]}{0pt}{}{#1}{#2}}

\begin{document}

\title{On the algebraic degree stability of vectorial Boolean functions when restricted to  affine subspaces
\thanks{The research of the first author is partly supported by the Norwegian Research Council; the research of the other two authors is supported by EPSRC, UK (EPSRC grant EP/W03378X/1).}}

\author{Claude Carlet\thanks{LAGA, Department of Mathematics, University of
Paris 8 (and Paris 13 and CNRS),
 Saint--Denis cedex 02, France, and University of Bergen, Norway.
E-mail: claude.carlet@gmail.com} \,  Serge Feukoua\thanks{
Department of Computer Science, University  of Loughborough, UK;
National Advanced School of Public Work, Cameroon.
E-mail:S.C.Feukoua-Jonzo@lboro.ac.uk; sergefeukoua@gmail.com}\, and
Ana S\u{a}l\u{a}gean\thanks{Departement of Computer Science,
University of Loughborough, UK. E-mail:A.M.Salagean@lboro.ac.uk}}
\maketitle \thispagestyle{empty}
\begin{abstract}
We study the behaviour of the algebraic degree of vectorial Boolean functions when their inputs are restricted to an affine subspace of their domain. Functions which maintain their degree on all subspaces of as high a codimension as possible are particularly interesting for cryptographic applications.

For functions which are power functions $x^d$ in their univariate representation, we fully characterize the exponents $d$ for which  the algebraic degree of the function stays unchanged when the input is restricted to spaces of codimension 1 or 2. For codimensions $k\ge 3$, we give a sufficient condition for the algebraic degree to stay unchanged. We apply these results to several known classes of power functions.
We define an optimality notion regarding the stability of the degree on subspaces, and determine a series of optimal functions, including the multiplicative inverse function and the quadratic APN functions.

We also give an explicit formula for counting the functions that keep their algebraic degree unchanged when restricted to hyperplanes.
\end{abstract}

 \noindent

\section{Introduction}
 Vectorial Boolean functions are   functions with $n$  input bits and $m$ output bits, where $n$ and $m$ are two positive integers
(see the  monograph \cite{Claude} for more details on vectorial Boolean  functions).
Such functions are used in symmetric cryptography in the design of block ciphers; they are usually called S-boxes in this context. They need to satisfy various conditions; in particular they need to have a high algebraic degree (denoted here $\deg(F)$) in order to avoid certain types of attacks such as higher-order differential attacks.
It is
important that the algebraic degree of these functions remains high even if the function is
restricted to an affine hyperplane or to an affine space of some small
codimension $k$ in order to avoid ``guess and determine attacks'', where the
attacker would make assumptions resulting in the fact that the input to the function is restricted to
a particular affine space.
Other attacks that exploit the behaviour (not necessarily the algebraic degree) of a function on a subspace of the original space  include the integral attack and the invariant subspace attack.

Note that the algebraic
degree of the restriction of a given function  to any affine
space is less than or equal to the algebraic degree of the global function and of course also less than or equal to the dimension of the affine space.
If the algebraic degree of a function $F$ restricted to an affine space $A$ of dimension $k$ is strictly lower than the algebraic degree of $F$, we will call $A$ a degree-drop space for $F$ ( note that when the degree of this restriction is at most zero, respectively, at most one, we retrieve the notion of $k$-normality, respectively, $k$-weak normality defined in  \cite{braeken2005normality}).
The functions which have no degree-drop space of dimension $k$, for $k$ as low as possible, are interesting for cryptographic purposes.
In \cite{CSA, CSA2}, the authors study this problem for Boolean functions (i.e.\ one output bit). They gave several characterizations, studied some particular classes (such as direct sums of monomials and the class of symmetric functions), explored the relationship to other parameters, and gave an explicit formula for the number of Boolean functions which maintain their degree on all hyperplanes.


In this paper, we
study this topic for vectorial Boolean functions. While some results from~\cite{CSA, CSA2} can be easily generalized from Boolean functions to vectorial Boolean functions by considering each coordinate function, (see Section~\ref{sec:basic}),  many other questions are more difficult and require different techniques,
for example by handling the univariate representation of the function.
Since we are interested in functions $F$ which have no degree-drop space of dimension $k$ for $k$ as low as possible, the best we can hope for is  $k=\deg(F)$ (as any function defined on a space of dimension $k$ has algebraic degree at most $k$);
 we say that such a function has \stab. While for Boolean functions this is impossible (it was proven in~\cite[Corollary~1]{CSA} that all Boolean functions $F$ have degree-drop spaces of dimension $\deg(F)$), vectorial Boolean functions with this property do exist.
 Namely, we will prove that the following functions have \stab: injective affine functions, quadratic APN function (for example a Gold APN function),  power function over $\F_{2^n}$ of the form $F(x) = x^{2^k-1}$ (including the multiplicative inverse function) and, more generally, $F(x) = x^{1+2^u+\ldots+2^{(j-1)u}}$, with $\gcd(u,n)=1$.

The related notion of $k$th order sum-freeness was introduced by the
first
author
in~\cite{Carlet2024,Carletadd}
, namely a function $F$ is $k$th order sum-free if for any affine space $A$ of dimension $k$, the sum of the values of $F$ over all elements of $A$ is non-zero.
Note that  the notion of $k$th order sum-freeness and  the notion of having no degree-drop space of dimension $k$ are different when $k\neq\deg(F)$, but when   $k=\deg(F)$ they coincide (and in this case it also coincides with having \stab). The functions with \stab\  listed above are therefore also $k$th order sum-free; this result was known for power functions over $\F_{2^n}$ of the form $F(x) = x^{2^k-1}$ (including the multiplicative inverse function), see~\cite[Proposition~1]{Carletadd}.
 Remark~\ref{rem:kth-order-sum-free-def} and~\ref{rem:kth-order-sum-free-results} provide further details regarding the connection to sum-freeness.



Of particular interest are functions which are power functions $x^d$ in their univariate representation (we have $n=m$ in this case). Such functions include the multiplicative inverse function $F(x)= x^{-1}$ (with $0^{-1} =0$ by convention), which can also be written as  $F(x)= x^{2^n-2}$; the S-box of the AES encryption algorithm is based on this function.

We fully characterize the exponents $d$ for which  
the function $x^d$ over $\F_{2^n}$
keeps its algebraic degree unchanged (i.e.\
has no degree-drop space) on spaces
of codimension 1 and 2 (see Theorems~\ref{thm:power-codim-1} and~\ref{thm:power-codim-2}). For codimensions $k\ge 3$, we give a sufficient condition,
see Theorem~\ref{thm:power-codim-k-u}.
Its proof relies on a result by Moore, \cite{Moore1896} and a result related to Moore determinants proven in~\cite{KshevetskiyGabidulin2005} in the context of MRD codes.
Essentially, Theorem~\ref{thm:power-codim-k-u} shows that it is sufficient that  a contiguous string of $k$ zeroes appears in the $n$-bit string of the binary representation of $d$  (considered cyclically), or after this string was decimated by an integer relatively prime to $n$. While this condition is sufficient, we show by means of counterexamples that it is not necessary.
This sufficient condition is further generalized in Theorem~\ref{thm:power-More-exponents} to require that the set of positions that contain zeroes in the binary representation of $d$ admits a subset of size $k$ which  is a ``Moore exponent set'' for $n$ (as defined in~\cite{BartoliZhou2020}).
We then apply these results to several classes of power functions, such as Gold APN functions, Kasami functions and  Welch functions.

For the multiplicative inverse function,
as the degree does decrease by at least $k-1$ when we restrict to spaces of codimension  $k\ge 2$, we go further and explore (Section~\ref{sec:inverse}) by how much it decreases. We show that it decreases by either $k-1$ or by $k$. Moreover, for $k=2$ we show that if $n$ is even, it always decreases by 1, while for $n$ odd, it decreases by~2 for a very small proportion (which we compute explicitly) of all the spaces. For $k=3$ we also describe the spaces for which the degree
decreases by 3.

Finally, in Section \ref{sec:cardofK1rnm}, we determine an explicit formula for the number of vectorial Boolean functions of a given degree $r$ in $n$ variables which keep their degree unchanged when restricted to any hyperplane. Our counting results use techniques developed in~\cite{SalMan17} and~\cite{SalFeru20}. A connection to functions which do not have ``fast points'' (i.e.\ functions $F$ whose discrete derivatives in any direction have the maximum possible degree, namely $\deg(F)-1$) was shown in~\cite{CSA2} for Boolean functions; we generalize this result
to vectorial Boolean functions. We then use this connection to also obtain an explicit formula for counting  the vectorial Boolean functions which do not have fast points, generalizing thus the result obtained in~\cite{SalMan17}. In~\cite{KolPol24} the authors exploit the fact that the sporadic Brinkmann-Leander-Edel-Pott function (which is the only known APN function
  that is not equivalent to either a monomial or a quadratic function) has $2^3-1$ fast points, which is a property that not many functions have; our counting results quantify how rare this property is: only about $\frac{1}{2^{103}}$ of the functions of degree 3 in 6 variables have this property.

\section{Preliminaries}\label{sec:prelim}
\subsection{Definitions and notation}\label{sec:def}
The finite field with two elements will be denoted as usual by $\mathbb{F}_2$. For every vector $v = (v_1, \ldots, v_n)\in\mathbb{F}_2^n$, the Hamming weight
$\w(v)$ of $v$ is the number of its non-zero coordinates i.e.\
the cardinality  of the support set
 $\{i\in \{1, \ldots, n\}:  v_i\neq 0\}$.  We shall denote by $|A|$ the cardinality of a set $A$.
 The all-zero vector of $\F_2^n$ will be denoted by $\mathbf{0}$.

A function $F:\mathbb{F}_2^n \rightarrow\mathbb{F}_2^m$ will be called a vectorial Boolean function, or simply an $(n,m)$-function. When $m=1$, $F$ is called a Boolean function.
Any $(n,m)$-function $F$ can be written as $F(x)=(f_1(x),...,f_m(x))$ with $f_i$ Boolean functions, called the coordinate functions of $F$.
Each $f_i$
can be uniquely represented in multivariate ANF (algebraic normal form), i.e.\ as a polynomial function  with coefficients in $\mathbb{F}_2$, in $n$ variables, and with degree at most 1 in each variable.
 The function $F$ itself can then be uniquely represented in multivariate ANF, i.e.\ as a polynomial of degree at most 1 in each variable, with coefficients in $\mathbb{F}_2^m$. In other words
 $F(x_1,\dots ,x_n)=\sum_{I\subseteq \{1,\dots ,n\}}a_I\prod_{i\in I}x_i$, where $a_I\in \F_2^m$
 (see~\cite{Claude}, for example).
 The {\em
algebraic degree} of $F$, denoted by $\deg(F)$, is
the degree of its ANF. Note that $\deg(F)=\max\{\deg(f_i):\,i=1,...,m\}$.  By convention, the algebraic degree of the identically zero function equals $-\infty$. A function $F$ of algebraic degree  at most 1
is called and {\em affine}
function; if an affine function also satisfies $F(\mathbf{0})=\mathbf{0}$, it is called {\em linear}.

%

Note that a second representation of an $(n,m)$-function
exists when $m=n$: endowing $\F_2^n$ with the structure of
the finite field with $2^n$ elements,  $\F_{2^n}$, any $(n,n)$-function  $F$ admits an
unique univariate polynomial  representation over
$\F_{2^n}$, of degree at most $2^n-1$ defined by
$F(x)=\sum_{j=0}^{2^n-1}c_jx^j$, with $c_i \in \F_{2^n}$. In this case, we
recover the algebraic degree of $F$ as  $\max\{\w(j): c_j\neq 0\}$ where for any integer $j$, by abuse of notation, we denote by $\w(j)$ the Hamming weight of the binary representation of $j$; in other words, if the binary representation of $j$ is $j=\sum_{t=0}^{n-1}j_t2^t$ with $j_t\in\{0,1\}$, then
$\w(j)=\sum_{t=0}^nj_t$.

The restriction of a function $F$ to an affine subspace $A$ of its domain will be denoted by $F_{\mid A}$. Note that, thanks to the affine identification between an affine space $A$ of codimension $k$ and $\mathbb{F}_2^{n-k}$, the function $F_{\mid A}$ can be viewed as a function in $n-k$ variables, which allows to consider its algebraic degree (which is independent of the choice of the affine identification) and we have then $\deg(F_{\mid A})\le \deg(F)$.

Recall that the number of vector subspaces of $\F_2^n$ of dimension $k$  is $\SqBinom{n}{k}_2$ (and it is also equal to the number of vector subspaces of codimension $k$), where the Gaussian $q$-binomial coefficients are defined for any $q>1$ as
\begin{equation}\label{eq:Gaussian-binom}
    \SqBinom{n}{k}_q = \frac{\prod_{i=n-k+1}^{n} (q^i-1)}{\prod_{i=1}^{k} (q^i-1)}.
\end{equation}

In this paper we are interested in functions which maintain their algebraic degree when restricted to spaces of a certain codimension.

\begin{definition}\label{deg-drop-def}
  Let $F$ be an $(n,m)$-function and $A$ an affine subspace of $\F_2^n$.
  If $\deg(F_{\mid A})< \deg(F)$, then we call $A$ a {\em degree-drop subspace} of $F$.
\end{definition}
Obviously, any subspace of a degree-drop space of a function $F$ is also a degree-drop space of $F$; if $F$ has no degree-drop space of dimension $k$, it has none of dimension $k+1$ either.

\begin{remark}
    Note that when $F$ has degree higher than one, if $A$ is a degree-drop space of dimension $k$ for $F$ and, moreover, $\deg(F_{\mid A})\le 0$ (or $\deg(F_{\mid A})\le 1$), then $F$ is called  $k$-normal (or $k$-weakly normal, respectively), see \cite{braeken2005normality}.
\end{remark}

Note also that $F_{|A}$ is a function in $\dim(A)$ variables, and therefore $\deg(F_{|A})\le \dim(A)$. This means that for any function $F$ and any $ k < \deg(F)$, any affine space of dimension $k$ is trivially a degree-drop space for $F$. The optimal functions, which have no other degree-drop spaces, are defined as follows:

\begin{definition}\label{def:stab}
    Let $F:\F_2^n \rightarrow \F_2^m$ be a vectorial Boolean function. We say that $F$ \emph{has \stab} if $\deg(F_{|A}) = \deg(F)$ for all affine spaces $A$ of $\F_2^n$ of dimension at least $\deg(F)$ (and it is sufficient for this that it happens for all affine spaces of dimension equal to $\deg(F)$).
\end{definition}

Note that in the case of Boolean functions (i.e.\ $m=1$) there are no functions which satisfy the condition above (it was proven in~\cite[Corollary~1]{CSA} that all Boolean functions $F$ have degree-drop spaces of dimension $\deg(F)$). 
For vectorial functions however, we shall see that such functions do exist. Namely, for degrees~1 and~2 we give such functions below in Propositions~\ref{prop:deg1} and~\ref{prop:deg2}. Further examples appear in Corollary~\ref{cor:monomial-k-stable}.

\begin{remark}\label{rem:kth-order-sum-free-def}
Note that the first author studied in \cite{Carlet2024,Carletadd} the dimensions $k$ of affine spaces $A$ over
which a given vectorial function $F$ (and in particular, the multiplicative  inverse function), sums to 0, that is, $\sum_{x\in A}F(x)=0$.
This is equivalent to the fact that
the restriction of $F$ to $A$ has degree strictly less than $k$.
Functions $F$ which do not admit such spaces for a given $k$ are said to be $kth$-order sum-free.
%
Note that for the  particular case of functions of degree $k$,
we have that $F$ is $k$th-order sum-free if and only if $F$ has no degree-drop space of dimension $k$, which is equivalent to $F$ having \stab. When $k\neq\deg(F)$ however, the notion of having no degree-drop space of dimension $k$ and the notion of being $k$th order sum-free are different. More precisely, when $k>\deg(F)$ the function $F$ cannot be $k$th-order sum-free, but there are functions $F$ which have no degree-drop space of dimension $k$; when  $k<\deg(F)$ the function $F$ will always have degree-drop spaces of dimension $k$, but there are functions $F$ which are $k$th-order sum-free.
Further connections between these notions are discussed in Remark~\ref{rem:kth-order-sum-free-results}.
\end{remark}

Recall that an affine automorphism of $\mathbb{F}_2^n$ is any mapping $\varphi$ of the form $\varphi(x) = Mx+a$, where $x$ and $a$ are viewed as column vectors $x^T=(x_1,\ldots, x_n)$, $a^T=(a_1,\ldots, a_n) \in \F_2^n$ and $M$ is an $n\times n$ invertible matrix over $\F_2$ (when the specific values of $M,a$ are relevant, we will also denote $\varphi$ as $\varphi_{M,a}$ or $\varphi_{M}$ if $a = \mathbf{0}$).
\begin{definition}\label{deffeqcomeq} Two $(n,m)$-functions $F$ and $G$ are said to be affinely equivalent if there  exist $\varphi$, an affine automorphism of $\mathbb{F}_2^n$ and $\phi$, an affine automorphism of $\mathbb{F}_2^m$,
such that $F=\phi\circ G\circ \varphi$ (or $F= G\circ \varphi$ when $F$ and $G$ are Boolean functions, meaning $\phi$ is the identity) where $\circ$ is the operation of
composition. We write then $F\sim G$.
\end{definition}

A parameter associated to a function is called an  affine invariant
if it is preserved by affine equivalence. The algebraic degree is an affine invariant. Moreover, in order to decide whether  a space $A$ is a degree-drop space for a function $F$ or degree $r$,
it suffices to examine the homogeneous function consisting of the monomials of degree $r$ of $F$; any monomials of lower degree can be ignored, so in effect it is sufficient to consider functions $F$ which are homogeneous. More precisely:
\begin{lemma}\label{lem:homog}
    Let $F$ be an $(n,m)$-function of degree $r$. Write $F=G+H$ such that $G$ is homogeneous of degree $r$ and $\deg(H)<r$. Then we have that an affine space $A$ is a degree-drop space for $F$ if and only if $A$ is a degree-drop space for $G$.
\end{lemma}
{\em Proof}.
We have $F_{|A} =G_{|A}+H_{|A}$, so $\deg(F_{|A})=r$ if and only if $\deg(G_{|A}+H_{|A})=r$; this in turn happens if and only if $\deg(G_{|A})=r$, since $\deg(H_{|A})\le \deg(H)<r$.
\hfill $\Box$\\

In view of the result above, we consider the natural extension of the affine equivalence $\sim$ between any $(n,m)$-functions to an equivalence $\sim_{r-1}$, which ignores monomials of degree $r-1$ or less.
Namely, for any fixed $r$ with $0\le r < n$, two $(n,m)$-functions $F$ and $G$ are equivalent under $\sim_{r-1}$ and we write  $F \sim_{r-1} G$, if and only if there is a function $H$ such that $F \sim G + H$ and $\deg(H)\le r-1$.

\subsection{Basic results}\label{sec:basic}

In this section, we will generalize a number of results from~\cite{CSA,CSA2} to vectorial Boolean functions. We will also settle the case of vectorial functions of degrees~1 and~2.

The property of having a degree-drop space is an affine invariant. More precisely:
 \begin{lemma}\label{lem:homomorphism}
Let $F$ and $G$ be two $(n,m)$-functions of degree $r$ such that  $F\sim_{r-1} G$ i.e.\
  $G = \phi \circ F \circ  \varphi + H$ for some affine automorphism $\phi$ of $\F_2^m$ and   $\varphi$ of $\F_2^n$
and some function $H$ with $\deg(H)\leq r-1$. Let $A$ be an affine space in $\F_2^n$.
Then $G_{\mid_A}\sim_{r-1} F_{\mid_{\varphi(A)}}$.
Therefore, an affine space $A$
is a degree-drop space for $G$ if and only if $\varphi(A)$ is a degree-drop space for $F$.
\end{lemma}
{\em Proof}.
The set $\varphi(A)$ is an affine space of the same dimension as $A$.
We have that $G_{\mid_A}=(\phi\circ F\circ \varphi + H)_{\mid_A}=((\phi\circ F + H\circ \varphi^{-1})\circ \varphi)_{\mid_A}$ is affine equivalent to $(\phi\circ F + H\circ \varphi^{-1})_{\mid_{\varphi(A)}}=(\phi\circ F)_{\mid_{\varphi(A)}} + (H\circ \varphi^{-1})_{\mid_{\varphi(A)}}$. Since  $(\phi\circ F)_{\mid_{\varphi(A)}}=\phi\circ (F_{\mid_{\varphi(A)}})$ is affine equivalent to $F_{\mid_{\varphi(A)}}$ and  $\deg((H\circ \varphi^{-1})_{\mid_{\varphi(A)}})\leq r-1$,
then we have $G_{\mid_A}\sim_{r-1} F_{\mid_{\varphi(A)}}$.
%
\hfill $\Box$\\

In order to determine whether a function has degree-drop affine spaces of a given dimension, it suffices to examine linear spaces only; applying~\cite[Lemma~5]{CSA2} to each coordinate function, we obtain:
\begin{lemma}\label{leaffivectreplace} Let $F$ be an $(n,m)$-function and  $A=a+E$ an affine space where $E$ is a vector space. Then $\deg(F_{\mid_A})=\deg(F)$ if and only if $\deg(F_{\mid_E})=\deg(F)$ (but note that when $\deg(F_{\mid_A})<\deg(F)$  we do not necessarily have $\deg(F_{\mid_A}) =\deg(F_{\mid_E})$).
\end{lemma}

We will now consider $(n,m)$-functions of degree one.
For Boolean functions, any function of degree one has a degree-drop hyperplane (see~\cite[Lemma~1(vii)]{CSA}). This is no longer the case for a vectorial function of degree~1. Noting that, when $F$ has degree one, a space $A$ is a degree-drop space for $F$ if and only if $F$ is constant on $A$, we can
easily determine more precisely the degree-drop spaces:
%
%
\begin{proposition}\label{prop:deg1}
    Let $F:\F_2^n \rightarrow \F_2^m$ be a linear function. The linear degree-drop spaces of $F$ are exactly those linear spaces $A$ with $A\subseteq \ker(F)$, where $\ker(F) = \{ x \in \F_2^n: F(x)=\mathbf{0} \}$ is the kernel of $F$.
    Therefore, an affine function $G$ has \stab\ if and only if $G$ is injective.
\end{proposition}

Recall that for a function $F:\F_2^n \rightarrow \F_2^m$  and $a\in\F_2^n\setminus \{\mathbf{0}\}$, the
discrete derivative of  $F$ in the direction $a$ is defined as
$D_aF(x)=F(x+a)+F(x)$.

A function $F:\F_2^n \rightarrow \F_2^n$ is called APN (almost perfect non-linear) if for any non-zero $a\in \F_2^n$ the derivative $D_aF$ is 2-to-1 (i.e.\ the preimage of any element of the codomain has cardinality at most 2).
Such functions have optimal resistance to differential cryptanalysis.
It is known (see e.g.\ \cite{NybergKnudsen1993}) that $F$
is APN  if and only if
the  restriction of $F$  to any space of dimension 2 is not affine. Therefore:
\begin{proposition}\label{prop:deg2}
  A quadratic  function $F:\F_2^n \rightarrow \F_2^n$ has no degree-drop space of dimension 2 (and therefore has \stab) if and only if it is an APN function.
\end{proposition}

\begin{remark}
    We could consider a stronger optimality condition than the one in  Definition~\ref{def:stab}, by having the additional requirement that $\deg(F_{|A}) = k$ for all spaces $A$ of dimension $k<\deg(F)$.
    Quadratic APN bijective functions, as well as linear injections (thanks to Proposition~\ref{prop:deg1}), have this property. We do not know whether any functions of degrees higher than two can satisfy this stronger optimality condition. However, we will see that the multiplicative inverse is ``close'' to satisfying this stronger condition, namely it satisfies $\deg(F_{|A}) \in\{k,k-1\}$ for all spaces $A$ of dimension $k<\deg(F)$, see Corollaries~\ref{cor:general-rest-deg-by-affine} and ~\ref{cor:general-rest-deg-by}
\end{remark}

We generalize a result from~\cite{CSA}  showing that the property that an $(n,m)$-Boolean function $F$ has a degree-drop space $A$ can be characterized using the indicator function of $A$. Recall that the indicator function, denoted by $1_A$, is the function $1_A: \F_2^n \rightarrow \F_2$ defined as $1_A(x)=1$ if $x\in A$ and $1_A(x)=0$ otherwise. We can then define the product $1_A F$ as the $(n,m)$-function defined as  $1_AF(x)=F(x)$ if $x\in A$ and $1_AF(x)=\mathbf{0}$ otherwise.

\begin{lemma}\label{add}
Let $F$ be a $(n,m)$-Boolean function and let $A$ be an affine subspace of $\mathbb{F}_2^n$, with indicator function $1_A$.
  We have: $$\deg (1_A F)=\deg(F_{\mid A})+\deg(1_A).$$
\end{lemma}
{\em Proof}. Let us set $F=(f_1,f_2,\ldots,f_m)$, where for all $i\in\{1,2,\ldots,m\}$, $f_i$ is a Boolean function. By observing that $F_{\mid A}=\left((f_1)_{\mid A},(f_2)_{\mid A},\ldots,(f_m)_{\mid A}\right)$ and that $(1_AF)=\left((1_Af_1),(1_Af_2),\ldots,(1_Af_m)\right)$), the proof is completed  by showing that the equalities $\deg (1_Af_i)=\deg((f_i)_{\mid A})+\deg(1_A)$ for all $i\in \{1,2,\ldots,m\}$ hold, which is true  thanks to \cite[Lemma~6]{CSA}.
%
%
\hfill $\Box$

Since the degree of $1_A$ is equal to the codimension of $A$, Lemma \ref{add} implies:
\begin{proposition}\label{p12}
Let $k$ be an integer and let $F$ be an $(n,m)$-function of algebraic degree $r$ such that $k\leq n$ and $1\leq r\leq n-k$. We have that $F$ has no degree-drop space of codimension $k$ if and only if, for all affine spaces $A$ of codimension $k$,
$$\deg(1_AF)=r+k.$$
\end{proposition}


Consider a homogeneous $(n,m)$-function $F=(f_1, \ldots,f_m)$. An affine space $A$ is a degree-drop space for $F$ if and only if $A$ is a degree-drop space for all the non-zero coordinate functions $f_i$. Consequently, a sufficient condition for $F$ to not have any degree-drop space of codimension $k$ is that at least one of its non-zero coordinate functions has no degree-drop space of codimension $k$ (for example, if one if its coordinate functions is a direct sum of $k+1$ monomials, see~\cite[Theorem~5]{CSA}). However, this condition is not necessary, as illustrated by the following example.

\begin{example}
     Consider $F(x)=(f_1(x),f_2(x),\ldots,f_n(x))$ of degree $n-1$ such that for all  $s\in\{1,2,...,n\}$,
$f_s(x)=x_1x_2\ldots x_{s-1}x_{s+1}\ldots x_n$. Clearly, each $f_s$ has degree-drop hyperplanes, namely all the hyperplanes defined by an affine equation in the variables $x_1, x_2, \ldots, x_{s-1}, x_{s+1}, \ldots, x_n$. Yet there is no hyperplane which is a degree-drop hyperplane for {\em all} the $f_s$; consequently  $F$ has no degree-drop hyperplane.
\end{example}

A number of further results from \cite{CSA} can be easily generalized from Boolean functions to vectorial Boolean functions:

\begin{lemma}\label{lem:x1-factor}
Let $F$  be a degree $r$ homogeneous $(n,m)$-function.
$F$ has a degree-drop hyperplane if and only if there exists $j\in \{1,\ldots,n\}$ and there exists a degree $r-1$ homogeneous $(n-1,m)$-function  $G$ 
such that $F(x)\sim_{r-1} x_jG(x_1,\ldots,x_{j-1},x_{j+1},\ldots, x_n)$.
\end{lemma}
{\em Proof.} 
Write $F(x)=(f_1(x),f_2(x),\ldots,f_m(x))$ and
let  $H$ be a degree-drop hyperplane  of $F$. This means  $H$ is a degree  drop hyperplane of all the non-zero
coordinates $f_i$  of $F$.
Thanks \cite[Theorem~2]{CSA} and its proof, $f_i=x_1g_i(x_2,\ldots, x_n)\circ \varphi +h_i(x_1,\ldots, x_n)$  where $g_i$ are homogeneous Boolean functions of degree $r-1$ in the $n-1$ variables $x_2,\cdots, x_n$,  $\deg(h_i)\leq r-1$ and $\varphi$ is the automorphism such that   $\varphi (H)=\{x:\,x_1=0\}$. Hence, by setting $G(x_2,\ldots, x_n)=(g_1(x_2,\ldots, x_n),\ldots,g_m(x_2,\ldots, x_n))$  and $H(x)=(h_1,\ldots,h_m)$, where  the  $(n,m)$-function $H$ is clearly of algebraic at most $r-1$,  we have then  $F(x)=x_1G_1(x_2,\ldots, x_n)\circ\varphi+H$ which
means $F(x)\sim_{r-1} x_1G(x_2,\ldots, x_n)$.
Note that  $x_1G(x_2,\ldots, x_n)\sim_{r-1} x_jG(x_1,\ldots,x_{j-1},x_{j+1},\ldots, x_n)$.

The converse is straightforward by considering the hyperplane defined by the equation $x_j=0$.
%
\hfill$\Box$

By applying~\cite[Theorem~3]{CSA} to each coordinate function  of a homogeneous $(n,m)$-function $F$ (by taking the same automorphism $\varphi$ for all the coordinate functions of $F$ which have degree $\deg(F)$ and then admit the same degree-drop affine space), we have the following result:
 \begin{lemma}\label{carcvect cod k}  Let $F$
 be a degree $r$ homogeneous $(n,m)$-function. $F$ has a degree-drop space of codimension $k$  if and only if $F\sim_{r-1}
G$ for some homogeneous function $G$  whose coordinate functions have an ANF in which each
monomial
contains at least one of the variables  $x_1,x_2,...,x_k$.

The affine homomorphism $\varphi$ such that $G=F\circ\varphi+H$ for some $H$ with $\deg(H)<r$, can be any affine transformation $\varphi$
mapping the vector space of equations $x_1=0$, $x_2=0$,..., $x_k=0$ to a degree-drop space of $F$ of codimension $k$.
\end{lemma}

\section{Power functions}\label{sec:power}
In this section, we work with $(n,n)$-functions and use their univariate representation. We are particularly interested in power functions  $F:\F_{2^n}\rightarrow \F_{2^n}$ defined as $F(x) = x^d$, with $1\le d \le 2^n-1$.
Since $L(x) = x^{2^j}$ is an invertible $\F_2$-linear function, we have that for any $j$ and any $d$ with $1\le d<2^n-1$ the functions $x^d$ and $x^{2^jd \bmod (2^n-1)}$  are affine equivalent.

We will characterize power functions with no degree drop space of codimension at most 2 and provide a sufficient condition for which  a power function has no degree drop space of codimension $k$ in general.  For the specific case of the inverse function $F(x) = x^{2^n-2}$ we analyze by how much the degree drops on spaces of codimensions~2 and~3.

Recall that the absolute trace function $\tr_n: \mathbb{F}_{2^n} \rightarrow \mathbb{F}_{2}$ is defined as $\tr_n(x)=x+x^{2}+x^{2^{2}}+...+x^{2^{n-1}}$. For any linear function  $F: \mathbb{F}_{2^n} \rightarrow \mathbb{F}_{2}$ there is an $a\in \F_{2^n}$ such that $F$ can be written as $F(x) = \tr_n(ax)$ (see, for example~\cite[Theorem~2.24]{LidNie94}). Any affine space $A$ of codimension $k$ can be expressed as the solution set of a system of $k$ independent affine equations $\tr_n(a_1x)+\epsilon_1=0, \ldots, \tr_n(a_1x)+\epsilon_1=0$ for some $\F_2$-linearly independent  $a_1, \ldots, a_k\in \F_{2^n}$ and some $\epsilon_1, \ldots, \epsilon_k\in \F_2$. Moreover, the indicator function of $A$ will be $1_A(x) = \prod_{i=1}^k (\tr_n(a_ix)+\epsilon_i +1)$. Therefore, Proposition~\ref{p12} and Lemma~\ref{leaffivectreplace} give:
\begin{corollary}\label{cor:indicator-trace}
    Let $F: \F_{2^n}\rightarrow \F_{2^n}$ and let $k\le n-\deg(F)$. We have that $F$ has no degree-drop space of codimension $k$ if and only if for all  $\F_2$-linearly independent  $a_1, \ldots, a_k\in \F_{2^n}$ we have $\deg\left(\tr_n(a_1x)\cdots \tr_n(a_kx) F(x)\right) = \deg(F)+k$.
\end{corollary}
While for an arbitrary $(n,m)$-function $F$ the degree of $F_{|A}$ can take any value lower than or equal to $\min(\deg(F), \dim(A))$, for the case when $F$ is a power function and $A$ is a linear space of codimension at most $\deg(F)$, the degree of $F_{|A}$ is lower bounded as follows:
\begin{proposition}\label{prop:deg-drop-at-most-k}
Let $F(x)=x^{d}$ be a power function over $\F_{2^n}$ of algebraic degree $r$, with $r<n$.
For every linear space $A$ of codimension $k$, with $1\le k \le r$, the degree of $F$ drops by at most $k$ over $A$ that is,
$r-k\leq \deg(F_{|A})$.
\end{proposition}
{\em Proof.}
Let $a_1,\ldots,a_k$ be linearly independent elements of $\F_{2^n}$ such that $A$ is defined by the equations $\tr_n(a_1x) =0, \ldots, \tr_n(a_kx) =0$. The indicator function of $A$ is therefore
\[1_A(x) = \prod_{i=1}^k (\tr_n(a_ix) +1)=1+G(x),\]
where $G(x)$ is a sum of functions of the form
$\tr_n(a_{i_1}x)\cdots \tr_n(a_{i_j}x)$ for some $j$ with $1\le j \le k$ and some distinct
$i_1, i_2, \ldots, i_j$ from  $\{1,\ldots,k\}$. Note that any product $\tr_n(a_{i_1}x)\cdots \tr_n(a_{i_j}x)$ has degree $j$ and there are no constant terms in its univariate representation.
 Therefore $G(x)$ has degree at most $k$ (and $k<n$) and has no constant term in its univariate representation; we can write $G(x)=\sum_{i=1}^{2^n-2}b_ix^i$ with $b_i\in\F_{2^n}$ (note that there are no terms $x^0$ and  $x^{2^n-1}$ in $G$).
Then, $1_{|A}(x)F(x)=x^d +\sum_{i=1}^{2^n-2}b_ix^{d+i}$. The term $x^d$ in the expression of $1_{|A}(x)F(x)$ cannot be cancelled by any other term (if it was cancelled, it would mean there was an $i$ with $x^d =  x^{d+i}$, which implies $i$ is a multiple of $2^n-1$, but that is impossible as $1\le i \le 2^n-2$). Therefore,  $\deg(1_{|A}F)\geq \deg(F)=r$ and
  by Lemma~\ref{add}, $\deg(F_{|A})=\deg(1_{|A}F)-k\geq r-k$.
%
 %
 \hfill$\Box$

 \begin{remark}
 An alternative proof of the lower bound in Proposition~\ref{prop:deg-drop-at-most-k} is as follows:
 Since  $x^d$ has algebraic degree $r$, the number $d’=2^n-1-d$ is such that $\w(d')=n-r$.
 The function $x^d x^{d’}$ equals $x^{2^n-1}$ and therefore we have $\sum_{x\in A} x^d x^{d’}=1\neq 0$. Let us decompose $F(x)=x^d$ over a basis $(\alpha_1,\ldots,\alpha_n)$ of $\F_{2^{n}}$: $F(x)=\sum_{i=1}^n f_i(x)\alpha_i$. Since $\sum_{i=1}^n \alpha_i (\sum_{x\in A} f_i(x) x^{d’}) \neq 0$, then there exists $i$ such that $\sum_{x\in A} f_i(x) x^{d’}\neq 0$. Then, there exists at least one element a such that \begin{align*}
     \tr_n(a\sum_{x\in A} f_i(x) x^{d’})&=\sum_{x\in A} f_i(x) \tr_n(ax^{d’})\nonumber\\&=\sum_{x\in F_2^n}1_A(x)f_i(x)\tr_n(ax^{d’})\neq 0.\end{align*} The Boolean function $\tr_n(ax^{d’})$ having an algebraic degree of at most $n-r$, the function $1_A(x)f_i(x)$ has then an algebraic degree of at least $n-(n-r)=r$ and  the proof is completed  by Lemma~\ref{add}.
 \end{remark}
\subsection{Characterization of power functions without degree-drop spaces}\label{sub-sec:power}

We start by looking at restrictions to spaces of codimension~1, then~2, followed by arbitrary codimensions $k$.
\begin{theorem}\label{thm:power-codim-1}
  Let $F:\F_{2^n}\rightarrow  \F_{2^n}$ defined as $F(x) = x^d$, with $1\le d \le 2^n-1$. Then $F$ has a degree-drop affine space of codimension 1  if and only if $d= 2^n-1$ (or equivalently, $\deg(F)=n$).
\end{theorem}
{\em Proof.} Note that if $d=2^n-1$ that is, $\deg(F)= n$, then $F$  all spaces of
codimension $1$ are degree-drop spaces for $F$ (because every vectorial function on a space of dimension $n-1$ has algebraic degree at most $n-1$). Now assume  $d\neq 2^n-1$ that is, $\deg(F)\le n-1$.
By Corollary~\ref{cor:indicator-trace},
it suffices to show that for all $a \in \F_{2^n}^*$ the product  $\tr_n(ax)x^{d}$  is of algebraic degree $\deg(F)+1$.  We have $\tr_n(ax)x^{d}=\sum_{t=0}^{n-1}a^{2^t}x^{2^t+d}$, where the exponents of $x$ are reduced modulo $2^n-1$, to the representatives $\{1,2, \ldots, 2^n-1\}$ (and not the representatives $\{0, 1,2, \ldots, 2^n-2\}$).
The terms in this sum cannot cancel each other, as for any two values $0\le t_1<t_2\le n-1$, if we assume, for a contradiction, that the exponents of  $x^{2^{t_1}+d}$ and $x^{2^{t_2}+d}$ are equal modulo $2^n-1$, it would mean that $2^{t_2}-2^{t_1}$ is a multiple of $2^n-1$, which is impossible.
Since $\deg(F)\le n-1$, in the binary representation of $d = \sum_{i=0}^n d_i 2^i$ we must have $d_j=0$ for at least one index $j$. We have $\w(2^j+d) = \w(d)+1$, therefore  The $a^{2^j}x^{2^j+d}$ is therefore of algebraic degree $\deg(F)+1$ and has a non-zero coefficient.
\hfill$\Box$


\begin{theorem}\label{thm:power-codim-2}
  Let $F:\F_{2^n}\rightarrow  \F_{2^n}$ be defined as $F(x) = x^d$, with $1\le d \le 2^n-1$ and $\deg(F)\le n-2$. Write the exponent $d$ in base 2, as $d = \sum_{i=0}^{n-1} d_i 2^i$ with $d_i\in\{0,1\}$. Let $u = \gcd(\{t_2-t_1: d_{t_1} =0, d_{t_2} =0, t_1\neq t_2\})$. We have that
   $F$ has no degree-drop affine space of codimension 2 if and only if $\gcd(u,n)=1$.
\end{theorem}
{\em Proof.}
By Corollary~\ref{cor:indicator-trace},
it suffices to show that for all   $a, b \in \F_{2^n}^*$, $a\neq b$, the product $\tr_n(ax)\tr_n(bx)x^{d}$ is of algebraic degree $\deg(F)+2$. We have
\begin{eqnarray*}
  \tr_n(ax)\tr_n(bx)x^{d} &=&  \sum_{t_1=0}^{n-1}\sum_{t_2=0}^{n-1}a^{2^{t_1}}b^{2^{t_2}}x^{2^{t_1}+2^{t_2}+d},
\end{eqnarray*}
 and again, the exponents of $x$ are reduced modulo $2^n-1$ to the representatives $\{1,2, \ldots, 2^n-1\}$ (note that for different pairs $\{t_1,t_2\}$ and $\{t’_1,t’_2\}$ we have $2^{t_1}+2^{t_2}\neq 2^{t'_1}+2^{t’_2}$ in the ring $\Bbb{Z}/(2^n-1)\Bbb{Z}$).
 Let $U_d = \{i\in \{0,1, \ldots, n-1\}: d_i=0\}$. Note that $x^{2^{t_1}+2^{t_2}+d}$ has degree $\w(d)+2$ if and only if  $t_1,t_2\in U_d$ and $t_1\neq t_2$. Therefore there is a function $G$ of degree at most $\deg(F)+1$ such that
\begin{align*}
  \tr_n(ax)\tr_n(bx) x^{d}&=  \sum_{\stackrel{t_1<t_2}{t_1,t_2\in U_d}}\left( a^{2^{t_1}}b^{2^{t_2}} + a^{2^{t_2}}b^{2^{t_1}}\right)x^{2^{t_1}+2^{t_2}+d}+G(x),
\end{align*}

where  the exponents $2^{t_1}+2^{t_2}+d$ are distinct for different $t_1,t_2$. It remains to determine what conditions on $d$ are necessary and sufficient to ensure that for all   $a, b \in \F_{2^n}^*$ with  $a\neq b$, there is at least one pair of integers 
$t_1,t_2\in U_d$, $t_1<t_2$ such that $a^{2^{t_1}}b^{2^{t_2}} + a^{2^{t_2}}b^{2^{t_1}}\neq 0$.

Let us consider the opposite, i.e.\ determine what conditions on $d$ are necessary and sufficient to ensure that there is at least one pair $(a,b)$ with   $a, b \in \F_{2^n}^*$ and  $a\neq b$ such that for all $t_1,t_2\in U_d$ with  $t_1<t_2$ we have $a^{2^{t_1}}b^{2^{t_2}} + a^{2^{t_2}}b^{2^{t_1}}= 0$.
Since $a$ and $b$ are non-zero, the equation $a^{2^{t_1}}b^{2^{t_2}} + a^{2^{t_2}}b^{2^{t_1}}= 0$ holds if and only if $\left(\frac{a}{b}\right)^{2^{t_2-t_1}-1} =1$.
Therefore we need that $\frac{a}{b} \in \F_{2^{t_2-t_1}} \setminus\F_{2}$ for all $0\le t_1<t_2\le n-1$ with $t_1,t_2\in U_d$.
This means $\frac{a}{b} \in \F_{2^{u}} \setminus\F_{2}$. Such a value for $\frac{a}{b}$ exists if and only if $\F_{2^{u}}\cap \F_{2^{n}} \neq \F_{2}$,
which in turn happens if and only if $\gcd(u,n)\neq 1$.
\hfill$\Box$

\begin{corollary}\label{cor:codim2}
  With the notation of Theorem~\ref{thm:power-codim-2}, if any of the following conditions is satisfied, then $F$ has  no degree-drop affine space of codimension 2:\\
  (i) There is an integer $t$ such that $0\le t \le n-1$ and $d_t=d_{(t+1)\bmod n}=0$.\\
  (ii) There are two  distinct integers $t_1, t_2\in \{0, 1, \ldots, n-1\}$ such that $d_{t_1} =d_{t_2}  = 0$ and $\gcd(t_2-t_1, n)=1$.\\
  (iii) There are three  distinct integers $t_1, t_2, t_3\in \{0, 1, \ldots, n-1\}$ such that $d_{t_1} =d_{t_2} =d_{t_3} = 0$ and $\gcd(t_2-t_1, t_3-t_2, n)=1$.\\
  (iv) $d<2^{n-2}$.\\
  (v) $\deg(F)\le \lfloor \frac{n-1}{2} \rfloor$.
\end{corollary}
{\em Proof.} \\
(i) One can check that the parameter $u$ in Theorem~\ref{thm:power-codim-2}, equals 1 and therefore, $\gcd(u,n)=1$.\\
(ii) $\gcd(t_2-t_1,  n)=1$ implies clearly $\gcd(u,n)=1$ in Theorem~\ref{thm:power-codim-2}.\\
(iii) $\gcd(t_2-t_1, t_3-t_2, n)=1$ implies clearly $\gcd(u,n)=1$ in Theorem~\ref{thm:power-codim-2}.\\
(iv) Since $d_{n-2}=d_{n-1}=0$, we can apply point (i).
\\
(v) If $\deg(F)\le \lfloor \frac{n-1}{2} \rfloor$, the weight of $d$ is at most $\lfloor \frac{n-1}{2} \rfloor$. One can verify that any sequence of $n$ bits of weight $\lfloor \frac{n-1}{2} \rfloor$ or less has at least two adjacent zeroes (considered circularly, i.e.\ the first and last bits are considered adjacent). We apply then point~(i).
\hfill$\Box$

We now give a sufficient condition for $x^d$ to have no degree-drop space of codimension $k$ for arbitrary $k$. We first need to recall the following result:
\begin{theorem}(\cite[Theorem III]{Moore1896})
\label{thm:Moore2}
    Let $q$ be a prime power. In $\F_{q}[x_1,\ldots,x_k]$ we have the following identity:
 \[ \displaystyle\prod_{(c_1,\ldots,c_k)\in P_k(\F_q)} (c_1x_1+\cdots +c_kx_k) = \begin{vmatrix}
        x_1 & \ldots & x_k \\
        x_1^{q} & \ldots & x_k^{q} \\
        \vdots &  & \vdots \\
        x_1^{q^{k-1}} & \ldots & x_k ^{q^{k-1}}
    \end{vmatrix}
    \]
    where $P_k(\F_q)$ denotes the set of non-zero vectors $(c_1,\ldots,c_k)\in \F_q^k$ such that the first non-zero coordinate is equal to $1$.
\end{theorem}
A simple consequence of this result is:
\begin{corollary}\label{cor:Moore-lin-indep}
Let $q$ be a prime power and let $a_1, \ldots,a_k\in \F_{q^m}$, with $k\le m$. Then the determinant
\[
\begin{vmatrix}
        a_1 & \ldots & a_k \\
        a_1^{q} & \ldots & a_k^{q} \\
        \vdots &  & \vdots \\
        a_1^{q^{k-1}} & \ldots & a_k^{q^{k-1}}
    \end{vmatrix}
\]
is non-zero if and only if $a_1, \ldots,a_k$ are $\F_{q}$-linearly independent.
\end{corollary}
We also recall the following result which is  not very well-known in the finite field theory (for example, it does not seem to appear in \cite{mullen2013handbook}). We proved this result, but later found out that it was proven  in a coding theory paper, namely \cite[Lemma~3]{KshevetskiyGabidulin2005}; since our proof is different, and uses an interesting result regarding irreducible polynomials, we included it in the Appendix.
\begin{proposition}(\cite[Lemma~3]{KshevetskiyGabidulin2005})
\label{prop:Lin-indep}Let $q$ be a prime power. Let $n$ and $u$ be two positive integers with $\gcd(n,u)=1$. For all  $a_1,a_2,\ldots,a_t\in \F_{q^n}$ with $t\leq n$, viewing $a_1,a_2,\ldots,a_t$  as elements of $\F_{q^{nu}}$,  we have that   $a_1,a_2,\ldots,a_t$ are $\F_q$-linearly independent if and only if  they are  $\F_{q^u}$-linearly independent.
\end{proposition}
Proposition~\ref{prop:Lin-indep} combined with Corollary~\ref{cor:Moore-lin-indep} yields the following result, presented in \cite{KshevetskiyGabidulin2005} in the context of MRD codes:
\begin{corollary}(\cite[Theorem~1]{KshevetskiyGabidulin2005})\label{cor:Moore-lin-indep-base-field}
Let $n$ and $u$ be two positive integers such that $\gcd(n,u)=1$. Let
$a_1, \ldots,a_k\in \F_{q^n}\subseteq \F_{q^{nu}}$, with $k\le n$.
Then the determinant
\[
\begin{vmatrix}
        a_1 & \ldots & a_k \\
        a_1^{q^u} & \ldots & a_k^{q^u} \\
        \vdots &  & \vdots \\
        a_1^{q^{u(k-1)}} & \ldots & a_k^{q^{u(k-1)}}
    \end{vmatrix}
\]
is non-zero if and only if $a_1, \ldots,a_k$ are $\F_{q}$-linearly independent.
\end{corollary}

We are now ready to give a sufficient condition for a power function $F(x)=x^d$ to not have degree-drop spaces of codimension $k$:
\begin{theorem}\label{thm:power-codim-k-u}
  Let $F:\F_{2^n}\rightarrow  \F_{2^n}$ defined as $F(x) = x^d$, with $1\le d \le 2^n-1$ and let $k$ be such that $k\le n-\deg(F)$. Write the exponent $d$ in base 2, as
  $d = \sum_{i=0}^{n-1} d_i 2^i$ with $d_i\in\{0,1\}$.
  If there
  are integers $t$ and $u$ such that $0\le t \le n-1$, $\gcd(u,n)=1$ and $d_t=d_{(t+u)\bmod n}=\cdots = d_{(t+u(k-1))\bmod n}=0$,
  then  $F$ has no degree-drop affine space of codimension $k$.
\end{theorem}
{\em Proof.}
%
By Corollary~\ref{cor:indicator-trace},
it suffices to show that for all $\F_2$-linearly independent $a_1, \ldots, a_k \in \F_{2^n}^*$ the product $x^{d}\prod_{i=1}^{k} \tr_n(a_ix)$ is of algebraic degree $\deg(F)+k$.  We have

 \[ x^{d}\prod_{i=1}^{k} \tr_n(a_ix) =  \sum_{t_1, \ldots, t_k\in\{0,\ldots,n-1\}}a_1^{2^{t_1}}\cdots a_k^{2^{t_k}}x^{2^{t_1}+\cdots+2^{t_k}+d}
\]
 and again, the exponents of $x$ are reduced modulo $2^n-1$ to the representatives $\{1,2, \ldots, 2^n-1\}$.
 Let $U_d = \{i\in \{0,1, \ldots, n-1\}: d_i=0\}$.
 Note that $x^{2^{t_1}+2^{t_2}+...+2^{t_k}+d}$ has degree $\w(d)+k$ if and only if  $t_1,t_2,...,t_k\in U_d$. Therefore there is a function $G$ of degree at most $\deg(F)+k-1$ such that the product $x^{d}\prod_{i=1}^{k} \tr_n(a_ix) $ equals
 \begin{align*}
  & \sum_{\stackrel{t_1<\cdots <t_k}{t_1,\ldots,t_k\in U_d}}\left( \sum_{\sigma\in S_k} a_{\sigma(1)}^{2^{t_1}}a_{\sigma(2)}^{2^{t_2}}\cdots a_{\sigma(k)}^{2^{t_k}}  \right)x^{2^{t_1}+\cdots+2^{t_k}+d}+G(x),
\end{align*}

where $S_k$ is the symmetric group consisting of all the permutations of $\{1,\ldots,k\}$.
Note that the exponents of $x^{2^{t_1}+\cdots+2^{t_k}+d}$ (reduced modulo $2^n-1$ as explained above) are distinct for  each distinct  tuple  $(t_1,\cdots ,t_k)$ such that  $t_1<\cdots <t_k$ (according to the uniqueness of the binary expansion of integers).
We need to show that at least one of them has a non-zero coefficient.

We examine the monomial obtained for $t_1=t, t_2=(t+u)\bmod n, \ldots, t_k = (t+u(k-1))\bmod n$ with $\gcd(n,u)=1$. It suffices to consider the case $t=0$
(otherwise we replace the function $x^d$ with the affine equivalent function $x^{2^{n-t}d}$).
The coefficient of $x^{1+2^u+2^{2u}+ \cdots +2^{u(k-1)}+d}$ is
\[
\sum_{\sigma \in S_k}  a_{\sigma(1)} a_{\sigma(2)}^{2^u} \cdots a_{\sigma(k)}^{2^{u(k-1)}},  \]
which is equal to the determinant
\[\begin{vmatrix}
        a_1 & \ldots & a_k \\
        a_1^{2^u} & \ldots & a_k^{2^u} \\
        \vdots &  & \vdots \\
        a_1^{2^{u(k-1)}} & \ldots & a_k ^{2^{u(k-1)}}
    \end{vmatrix}.
\]
By Corollary~\ref{cor:Moore-lin-indep-base-field}, this determinant is non-zero, as $a_1, \ldots, a_k$ are   $\F_2$-linearly independent.
Therefore  the coefficient of {\small$x^{1+2^u+2^{2u}+ \cdots +2^{u(k-1)}+d}$}  is non-zero. Since its degree is $\deg(F)+k$, this ends the proof.
\hfill$\Box$

The condition in Theorem~\ref{thm:power-codim-k-u} is sufficient, but not necessary for not having degree-drop spaces of codimension $k$, as shown by the next example.
\begin{example}\label{ex:counterexample}
    For $k=2$, consider $n=70$ and the power function $x^d$ where $d=\sum_{i=0}^{69}d_i2^i$, with $d_i=0$ when $i\in\{0,6,21\}$ and $d_i=1$ otherwise. Applying Theorem~\ref{thm:power-codim-2}, we see that $F$ has no degree-drop space of codimension~2. However, putting $t_1=0, t_2=6, t_3=21$, we note that $\gcd(t_i-t_j,n)>1$ for all distinct $i,j$ in $\{1,2,3\}$, and therefore the sufficient condition of Theorem~\ref{thm:power-codim-k-u} is not satisfied.

    For $k=3$, consider $n=8$ and $F(x) = x^{39}$. One can easily check that this function does not satisfy the sufficient condition in Theorem~\ref{thm:power-codim-k-u}. However, a computer calculation shows that $F$ does not have degree-drop spaces of codimension~3.
\end{example}
Theorem~\ref{thm:power-codim-k-u} allows us to construct a family of functions which have \stab.
\begin{corollary}\label{cor:monomial-k-stable}
     Let $F:\F_{2^n}\rightarrow  \F_{2^n}$ be a power function of algebraic degree $j$ defined as $F(x) = x^d$, with $d=1+2^u+\ldots +2^{(j-1)u}$, where $u$ and $j$ are integers such that $\gcd(u,n)=1$ and $1<j<n$. Then $F$ has no degree-drop space of codimension $n-j$ and therefore $F$ has \stab. (In particular for $u=1$ we obtain that $F(x)=x^{2^j-1}$ has \stab; when $u=1$ and $j=n-1$, $F(x)=x^{2^{n-1}-1}$ is affine equivalent to the multiplicative inverse function.).
\end{corollary}
{\em Proof.}
We apply Theorem~\ref{thm:power-codim-k-u} with $t=ju$ and $k = n-j$, and  $u$ in  Theorem~\ref{thm:power-codim-k-u} being the same as the $u$ here.  By Definition~\ref{def:stab}, $F$ has \stab.
\hfill$\Box$
\begin{remark}\label{rem:kth-order-sum-free-results}
As mentioned in Remark~\ref{rem:kth-order-sum-free-def}, for a function $F$ of degree $k$, we have that $F$ is $k$th order sum-free if and only if $F$ has \stab. It was proven in~\cite{Carlet2024} that $F(x)=x^{2^j-1}$ is $j$th order sum-free, which means $F$ has \stab. This result is retrieved as the particular case of $u=1$ in Corollary~\ref{cor:monomial-k-stable},
which  also provides a new, more general  class of $j$th order sum-free functions.
\end{remark}
\begin{corollary}\label{cor:monomial-k}
   Let $F:\F_{2^n}\rightarrow  \F_{2^n}$ defined as $F(x) = x^d$, with $1\le d \le 2^n-1$ and let $k$ be such that $k\le n-\deg(F)$.
  If any of the following conditions is satisfied, then $F$ has  no degree-drop affine space of codimension k:\\
  (i) $d<2^{n-k}$.\\
  (ii) $\deg(F)\le \lfloor \frac{n-1}{k} \rfloor$.\\
  then $F$ has \stab.
\end{corollary}
{\em Proof.} Use Theorem~\ref{thm:power-codim-k-u}  and apply the same type of argument as in the proof of Corollary~\ref{cor:codim2}(iv, v).
\hfill$\Box$

Next, we apply Theorem~\ref{thm:power-codim-k-u} to several well-known classes of functions.
\begin{example}
    Consider the quadratic Gold function $F(x) = x^{1+2^j}$ over $\F_{2^n}$, which is APN if and only if $\gcd(n,j)=1$. When $\gcd(n,j)=1$, the conditions of Theorem~\ref{thm:power-codim-k-u} are satisfied for $t=2j$, $u=j$ and $k=n-2$. Therefore, $F$ has no degree-drop space of codimension $n-2$ (i.e.\ dimension~2). The same result follows from Proposition~\ref{prop:deg2}.
\end{example}
\begin{example}
 Another class of APN functions are  the Welch functions~\cite{Dobbertin1999}. A Welch function is defined over $\F_{2^n}$ when $n$ is odd, as $F(x)=x^{2^j+3}$ with $j = \frac{n-1}{2}$. It has algebraic degree~3.
 Since  $\gcd(n,j)=1$,  the conditions of Theorem~\ref{thm:power-codim-k-u} are satisfied for $t=2j$, $u=j$ and $k=n-4$. Therefore, $F$ has no degree-drop space of codimension $n-4$ (i.e.\ dimension~4), which is close to the optimal value, namely $n-\deg(F) = n-3$.

    \end{example}

\begin{example}
 Let us consider the case of
 the Kasami function
 defined by $F(x)=x^d$
 with $d=2^{2i}-2^{i}+1$.
Note that when $\gcd(i,n)=1$,
this class of power functions
is included in the class of APN functions.
Examining the binary representation of the exponent $d=2^{2i}-2^{i}+1$ we see that the algebraic degree is $i+1$, and there are two blocks of successive zeroes, one of indices $1,\dots,i-1$, of length $i-1$, and one of indices  $2i,\dots,n-1$, of length $n-2i$. Therefore, by Theorem~\ref{thm:power-codim-k-u},
the Kasami function  has no degree drop space of codimension $k= \max(i-1, n-2i)$.

    \end{example}
The property in Corollary~\ref{cor:Moore-lin-indep} inspired the following definition:
\begin{definition}(\cite{BartoliZhou2020})
    Let $q$ be a prime power and  let $a_1, \ldots, a_k\in \F_{q^n}$.
  A set   $\{i_0, \ldots, i_{k-1}\}$ of non-negative integers is called a {\em Moore exponent set} for $q$ and $n$ if the following property holds: the determinant
   \[
  \begin{vmatrix}
        a_1^{q^{i_0}} & \ldots & a_k^{q^{i_0}} \\
        a_1^{q^{i_1}} & \ldots & a_k^{q^{i_1}} \\
        \vdots &  & \vdots \\
        a_1^{q^{i_{k-1}}} & \ldots & a_k^{q^{i_{k-1}}}
    \end{vmatrix}
  \]
  is non-zero if and only if $a_1, \ldots, a_k$ are $\F_q$-linearly independent.
\end{definition}
A set  $\{i_0, \ldots, i_{k-1}\}$ is a Moore exponent set for $q$ and $n$  if and only if $\{0, (i_1-i_0)\bmod n \ldots, (i_{k-1}-i_0)\bmod n\}$ is a Moore exponent set for $q$ and $n$ (see~\cite{BartoliZhou2020}).

From Theorem~\ref{thm:Moore2} we see that $\{0,1,\ldots,k-1\}$ with $k\le n$ is a Moore exponent set for any prime power $q$ and any $n$; more generally, from Corollary~\ref{cor:Moore-lin-indep-base-field},  when
$\gcd(n,u)=1$ and $k\le n$, the set $\{0,u,\ldots,(k-1)u\}$ is a Moore exponent set for $q$ and $n$.

 We can generalize Theorem~\ref{thm:power-codim-k-u} to Moore exponent sets as follows:
 \begin{theorem}\label{thm:power-More-exponents}
  Let $F:\F_{2^n}\rightarrow  \F_{2^n}$ defined as $F(x) = x^d$, with $1\le d \le 2^n-1$ and let $k$ be such that $k\le n-\deg(F)$. Write the exponent $d$ in base 2, as
  $d = \sum_{i=0}^{n-1} d_i 2^i$ with $d_i\in\{0,1\}$.
  If there is a subset of cardinality $k$ of $\{i: d_i=0\}$ which is a Moore exponent set for $2$ and $n$,
  then  $F$ has no degree-drop affine space of codimension $k$. The converse does not hold.
\end{theorem}
 {\em Proof.}
 Let $\{i_0, \ldots, i_{k-1}\}$ be a Moore exponent set which is a subset of $\{i: d_i=0\}$. The proof is similar to the proof of Theorem~\ref{thm:power-codim-k-u}, but this time we examine
the coefficient of $x^{2^{i_0}+ \cdots +2^{i_{k-1}}+d}$, namely
\[
\sum_{\sigma \in S_k}  a_{\sigma(1)}^{2^{i_0}}  \cdots a_{\sigma(k)}^{2^{i_{k-1}}},  \]
which is equal to the determinant
\begin{equation}\label{eq:Moore-exponent-set-det}
    \begin{vmatrix}
        a_1^{2^{i_0}} & \ldots & a_k^{2^{i_0}} \\
        a_1^{2^{i_1}} & \ldots & a_k^{2^{i_0}} \\
        \vdots &  & \vdots \\
        a_1^{2^{i_{k-1}}} & \ldots & a_k^{2^{i_{k-1}}}
    \end{vmatrix},
\end{equation}
which is non-zero as $\{i_0, \ldots, i_{k-1}\}$ is a Moore exponent set.

For the converse, note that the fact that $F$ has no degree-drop space of codimension $k$ implies that for all $\F_2$-linearly independent $a_1, \ldots, a_k\in \F_{2^n}$  there is a subset $\{i_0, \ldots, i_{k-1}\}$ of $\{i: d_i=0\}$ such that the determinant~\eqref{eq:Moore-exponent-set-det} is non-zero. However, it is not necessary that the same set $\{i_0, \ldots, i_{k-1}\}$ works for all $\F_2$-linearly independent $a_1, \ldots, a_k\in \F_{2^n}$ (if there is a set that works for all, then that set would be a Moore exponent set). Indeed we can give a counterexample.
Consider $k=2$. It is easy to check that $\{i_0, i_1\}$ is a Moore exponent set for $q$ and $n$ if and only if
$\gcd(i_1-i_0,n)=1$. In other words, for $k=2$, Theorem~\ref{thm:power-More-exponents} coincides with Theorem~\ref{thm:power-codim-k-u}.
The first power function in Example~\ref{ex:counterexample},
which is a counterexample for the converse of Theorem~\ref{thm:power-codim-k-u}, is therefore also a counterexample for the converse of Theorem~\ref{thm:power-More-exponents}.
\hfill$\Box$

Note that Theorem~\ref{thm:power-codim-k-u} becomes a particular case of Theorem~\ref{thm:power-More-exponents}, for the Moore exponent set $\{t,t+u,\ldots,t+(k-1)u\}$, with $\gcd(n,u)=1$.

\subsection{
The multiplicative inverse function}\label{sec:inverse}
Consider the multiplicative inverse
function $I:\F_{2^n}\rightarrow \F_{2^n}$ defined as $I(x) = x^{-1}$, with the convention $0^{-1}=0$. Note that this function can also be expressed as a power function,
$I(x)=x^{2^n-2}$, and its algebraic degree is $n-1$. It is also affine equivalent to the function $x^{2^{n-1}-1}$ so, applying Theorem~\ref{thm:power-codim-1}
 or Corollary~\ref{cor:monomial-k}, we obtain
\begin{corollary}
    The multiplicative inverse function $I:\F_{2^n}\rightarrow \F_{2^n}$ defined as $I(x) = x^{-1}$ has no degree-drop hyperplanes and has \stab.
\end{corollary}
 We can easily observe that
any space of codimension 2 or more is a degree-drop space for the multiplicative inverse function (because every vectorial function on a space of dimension $n-2$ has algebraic degree at most $n-2$).

In \cite[Theorem 1]{Carletadd}, it was proven that for any affine space $A$ which is not a linear space, the sum $\sum_{x\in A} I(x)$ is non-zero, which is equivalent to the fact that the degree of $I_{|A}$ equals the dimension of $A$ (see Remark~\ref{rem:kth-order-sum-free-def}). We obtain therefore:
\begin{corollary}  (see~\cite[Theorem 1]{Carletadd})
    \label{cor:general-rest-deg-by-affine}
    Let $I$  the multiplicative inverse function over $\F_{2^n}$. Let $A$ be an affine subspace of $\F_{2^n}$ of codimension $k$ with $2\leq k<n$ such that $A$ is not a linear space. Then the degree of $I$ drops by   $k-1$
    over $A$ (i.e.\  $\deg(I_{|A})=n-k$).
\end{corollary}

For the case when $A$ is a linear space, Proposition~\ref{prop:deg-drop-at-most-k}, combined with the obvious fact that $\deg(F_{|A})\le \dim(A)$ gives:
\begin{corollary}\label{cor:general-rest-deg-by} Let $A$ be a linear subspace of $\F_{2^n}$ of codimension $k$ with $2\leq k<n$ and $I$  the multiplicative inverse function. Then the degree of $I$ drops by   $k-1$ or $k$ over $A$ (i.e. $\deg(I_{|A})\in\{n-k,n-1-k\}$).
\end{corollary}
 \begin{remark}
 It was conjectured  \cite{Carletadd}  that for every \( n \) and every \( k \), \( 3 \leq k\leq  n-3 \), the inverse function is not \( k \)th-order sum-free. This conjecture has been proven in each of the following cases (but remains open for the other cases):
(1) \( k \) has a non-trivial common divisor with \( n \) (see~\cite{Carletadd}),
(2) \( n \) is even,
(3) \( k \) or \( n-k \) is less than \( \frac{n}{11} \) (for (2) and (3) see~\cite{CarletHou2024}). Other conditions for which this conjecture is true can be found in~\cite{HouZhao2025,HouZhao2-2025, EbelingetAl2024}.  In our context, this conjecture is equivalent to: ``for every \( k \), the degree of \( I \) drops by at least \( k \) on at least one space \( A \) of codimension \( k \)''.
Corollary~\ref{cor:general-rest-deg-by} provides a refinement of this conjecture by proving that it never drops by more than $k$.
 \end{remark}

When $k=2$ or $3$, the following results provide  more precision  regarding the spaces on which the degree drops by $k-1$ and the ones on which it drops by  $k$.
\begin{theorem}\label{thm:inv-codim-2}
Let
$I$ be the multiplicative inverse function on $\F_{2^n}$.
We have:
\begin{itemize}
\item For $n$ odd, the algebraic degree  of $I$ drops by 1 on all spaces $A$ of codimension~$2$ (i.e.\ $\deg(I_{|A})=n-2$).
 \item For $n$ even, the algebraic degree  of $I$ drops by  2   on all  linear spaces $A$ of codimension $2$
 defined by the equations $\tr_n(ax)=0$  and $\tr_n(cax)=0$ with $a,c\in \F_{2^n}^*$ such that $c^2+c+1=0$
(and there are $\frac{2^n-1}{3}$ such spaces, representing a proportion of $\frac{1}{2^{n-1}-1}$ of all the linear spaces of codimension~2); it drops by 1 on all the other affine spaces $A$ of codimension $2$.
 \end{itemize}  
\end{theorem}
{\em Proof.}
%
According to  Lemma~\ref{add}, $\deg(I_{\mid A}) = \deg(1_A I)-2$. Let the affine space $A$ of codimension~2 be defined by the equations $\tr_n(ax)+\gamma_1=0$  and $\tr_n(bx)+\gamma_2=0$ where
 $a,b\in\F_{2^n}^*$, $a\neq b$ and    $\gamma_1,\gamma_2\in \F_2$.
 The indicator function of  $A$ is  $1_A(x) = (\tr_n(ax)+\epsilon_1)(\tr_n(bx)+\epsilon_2)$ where
$\epsilon_1=\gamma_1+1$ and $\epsilon_2=\gamma_2+1$.

We have
\begin{align*}
&(\tr_n(ax)+\epsilon_1)(\tr_n(bx)+\epsilon_2)I(x) = \sum_{i,j=0}^{n-1} a^{2^i}b^{2^j} x^{2^n+2^i+2^j-2}+\\&\epsilon_1
\sum_{j=0}^{n-1}b^{2^j} x^{2^n+2^j-2} +
\epsilon_2
\sum_{i=0}^{n-1}a^{2^i} x^{2^n+2^i-2}+\epsilon_1\epsilon_2I(x).
\end{align*}
(the exponents of $x$ are reduced modulo $2^n-1$, to the representatives $\{1,2, \ldots, 2^n-1\}$ and not the representatives $\{0, 1,2, \ldots, 2^n-2\}$).
The coefficient of $x^{2^n-1}$ in $(\tr_n(ax)+\epsilon_1)(\tr_n(bx)+\epsilon_2)I(x)$ equals
\[
a^{2^{n-1}}b^{2^{n-1}}+\epsilon_1b+\epsilon_2a
\]

We can assume, without loss of generality,
that $\epsilon_1\le \epsilon_2$ (as integers in $\{0,1\}$).
When $\epsilon_1 = \epsilon_2 =0$, the coefficient of $x^{2^n-1}$ is  $a^{2^{n-1}}b^{2^{n-1}}$, which is non-zero, so the algebraic degree of $1_A I$ is $n$. When $\epsilon_1 =0, \epsilon_2 =1$, the coefficient of $x^{2^n-1}$ is $(ab +a^2)^{2^{n-1}}$, which is non-zero, so, again, the algebraic degree of $1_A I$ is $n$. (Alternatively, these two cases follow from Corollary~\ref{cor:general-rest-deg-by-affine}).

When $\epsilon_1 = \epsilon_2 =1$, the coefficient of $x^{2^n-1}$ is $a^{2^{n-1}}b^{2^{n-1}} +a+b=(ab+a^2+b^2)^{2^{n-1}}$.
We have that $ab+a^2+b^2=0$  if and
only if $c^2+c+1=0$, with $c=b/a$. We know that $c^2+c+1=0$
is the equation satisfied by the elements of $\Bbb{F}_4\setminus \F_2$. So the
necessary and sufficient condition for the equation $c^2+c+1=0$ to have solutions in $\F_{2^n}$
is that $\Bbb{F}_4\subseteq \Bbb{F}_{2^n}$
that is, $n$ is even. This means that when $n$ is odd, the degree
of the function $I$ will decrease by only $1$ over all spaces of
codimension $2$. When $n$ is even, the equation  $c^2+c+1=0$ has exactly two solutions in $\F_{2^n}$, namely the elements of $\F_4\setminus\F_2$. This means that there are $2(2^n-1)$ solutions $(a,b)$ for the equation $ab+a^2+b^2=0$. For each vector space $V$ of codimension~2 there are $(2^2-1)(2^2-2) = 6$ ways to pick an (ordered) tuple of values $(a,b)$ such that $V$ is defined by the system of two equations $\tr_n(ax)=0$, $\tr_n(bx)=0$. The number of spaces $A$ on which the degree drops by at least~2 is therefore $\frac{2(2^n-1)}{6} = \frac{2^n-1}{3}$. Recalling that the total number of vector spaces of codimension~2 is $\SqBinom{n}{2}_2$, we can verify that $\frac{2^n-1}{3\SqBinom{n}{2}_2}=\frac{1}{2^{n-1}-1}$.

 In the case that the coefficient of $x^{2^n-1}$ vanishes (when $n$ is even), we have $\deg(I_{|A})\le n-3$, and
thanks to Proposition~\ref{prop:deg-drop-at-most-k}, $\deg(I_{|A})\ge n-3$, and therefore $\deg(I_{|A})= n-3$.
\hfill$\Box$

For the case of affine spaces of codimension 3, we have:
\begin{theorem}\label{Theo:inv-codim3} Let $A$ be a codimension $3$ affine space of $\F_{2^n}$  defined by the equations
$\tr_n(a_1x)+\gamma_1=0$, $\tr_n(a_2x)+\gamma_2=0$  and $\tr_n(a_3x)+\gamma_3=0$
where  $a_1,\, a_2$ and $a_3$ are three $\F_2$-linearly independent  elements of $\F_{2^n}$ and    $\gamma_1,\gamma_2,\gamma_3\in \F_2$. The algebraic degree of the multiplicative inverse function  $I$
\begin{itemize}
\item drops by  3
on all spaces $A$ which are linear (i.e.\  $(\gamma_1,\gamma_2,\gamma_3)= (0,0,0)$) and
\begin{equation}\label{eq:indicator-3}
  d_1d_2(1+d_1+d_2)+d_1^2+d_2^2+d_1^2d_2^2+(1+d_1^4+d_2^4)=0,
\end{equation}
where $d_1=a_2/a_1$ and
 $d_2=a_3/a_1$.
 \item drops by 2 on all other affine spaces of codimension~3.
\end{itemize}
%
%
\end{theorem}
{\em Proof.} Since $I$ is of algebraic degree $n-1$,  any space of codimension at least $2$ is a degree-drop space of $I$.
According to Lemma \ref{add},
$\deg(I_{\mid A}) = \deg(1_A I)-3$. The affine space $A$ of codimension 3 has the indicator function $1_A(x) = \prod_{s=1}^3(\tr_n(a_sx)+\epsilon_s)$
where $\epsilon_i=\gamma_i+1$ for all $i\in\{1,2,3\}$.
We have
 \begin{align*}
& I(x)\prod_{s=1}^3(\tr_n(a_sx)+\epsilon_s)=\sum_{i_1,i_2,i_3=0}^{n-1} a_{1}^{2^{i_1}}a_{2}^{2^{i_2}}a_{3}^{2^{i_3}}x^{2^n+2^{i_1}+2^{i_2}+2^{i_3}-2}+\\&\sum_{s=1}^3\epsilon_s\sum_{\substack{j< k\\j, k\in\{1,2,3\}\setminus\{s\}}}\sum_{i_j,i_k=0}^{n-1} a_{j}^{2^{i_j}}a_{k}^{2^{i_k}} x^{2^n+2^{i_j}+2^{i_k}-2} +\\&
\sum_{s=1}^3\sum_{\substack{j< k\\j, k\in\{1,2,3\}\setminus\{s\}}}\epsilon_{j}\epsilon_{k}\sum_{i_s=0}^{n-1} a_{s}^{2^{i_s}}x^{2^n+2^{i_s}-2}
%
+\epsilon_1\epsilon_2\epsilon_3I(x).
 \end{align*}

(the exponents of $x$ are reduced modulo $2^n-1$, to the representatives $\{1,2, \ldots, 2^n-1\}$ and not the representatives $\{0, 1,2, \ldots, 2^n-2\}$).
The coefficient of $x^{2^n-1}$ in $I(x)\prod_{s=1}^3(\tr_n(a_sx)+\epsilon_s)$ equals
\begin{eqnarray*}
& & a_1^{2^{n-1}}a_2^{2^{n-2}}a_3^{2^{n-2}} + a_1^{2^{n-2}}a_2^{2^{n-1}}a_3^{2^{n-2}} + a_1^{2^{n-2}}a_2^{2^{n-2}}a_3^{2^{n-1}} \nonumber \\
&+& \epsilon_1 a_2^{2^{n-1}}a_3^{2^{n-1}} + \epsilon_2 a_1^{2^{n-1}}a_3^{2^{n-1}} + \epsilon_3 a_1^{2^{n-1}}a_2^{2^{n-1}} \nonumber \\
&+& \epsilon_1 \epsilon_2 a_3 + \epsilon_1 \epsilon_3 a_2 + \epsilon_2 \epsilon_3 a_1 \nonumber \\
&=& a_1^{2^{n-1}}a_2^{2^{n-2}}a_3^{2^{n-2}} + a_1^{2^{n-2}}a_2^{2^{n-1}}a_3^{2^{n-2}} + a_1^{2^{n-2}}a_2^{2^{n-2}}a_3^{2^{n-1}} \nonumber \\
&+& \epsilon_1 a_2^{2^{n-1}}a_3^{2^{n-1}} + \epsilon_2 a_1^{2^{n-1}}a_3^{2^{n-1}} + \epsilon_3 a_1^{2^{n-1}}a_2^{2^{n-1}} \nonumber \\
&+& \epsilon_1 \epsilon_2 a_3^{2^{n}} + \epsilon_1 \epsilon_3 a_2^{2^{n}} + \epsilon_2 \epsilon_3 a_1^{2^{n}}
\end{eqnarray*}
%

We can assume without loss of generality that $\epsilon_1\le \epsilon_2\le \epsilon_3$ (as integers in $\{0,1\}$). When $\epsilon_1=\epsilon_2=\epsilon_3=0$, the coefficient of $x^{2^n-1}$ equals
$a_1^{2^{n-1}}a_2^{2^{n-2}}a_3^{2^{n-2}}+a_1^{2^{n-2}}a_2^{2^{n-1}}a_3^{2^{n-2}}+
a_1^{2^{n-2}}a_2^{2^{n-2}}a_3^{2^{n-1}}=(a_1a_2a_3(a_1+a_2+a_3))^{2^{n-2}}\neq 0$. For
$\epsilon_1=\epsilon_2=0 $ and $\epsilon_3=1$,
the coefficient of $x^{2^n-1}$ equals
  $(a_1a_2a_3(a_1+a_2+a_3)+a_1^2a_2^2)^{2^{n-2}}=(a_1a_2(a_3+a_1)(a_3+a_2))^{2^{n-2}}$; for $\epsilon_2=\epsilon_3=1 $ and $\epsilon_1=0$
it equals
$(a_1a_2a_3(a_1+a_2+a_3)+a_1^2a_2^2+a_1^2a_3^2+a_1^4)^{2^{n-2}}=(a_1(a_1+a_2+a_3)(a_1+a_3)(a_1+a_2))^{2^{n-2}}$,
which are all non-zero. (Alternatively, these three cases follow from Corollary~\ref{cor:general-rest-deg-by-affine}).

In the case $\epsilon_1=\epsilon_2=\epsilon_3=1$, the
 coefficient of  $x^{2^n-1}$ equals $A(a_1,a_2,a_3)^{2^{n-2}}$,
 where $A(a_1,a_2,a_3)=a_1a_2a_3(a_1+a_2+a_3)+a_1^2a_2^2+a_1^2a_3^2+a_2^2a_3^2+a_1^4+a_2^4+a_3^4$, which, taking into account that all $a_i$ are non-zero, can be rewritten as
 $A(a_1,a_2,a_3)=a_1^4\left(\frac{a_2}{a_1}\frac{a_3}{a_1}(1+\frac{a_2}{a_1}+\frac{a_3}{a_1})+\frac{a_2^2}{a_1^2}+\frac{a_3^2}{a_1^2}+\frac{a_2^2a_3^2}{a_1^4}+1+\frac{a_2^4}{a_1^4}+\frac{a_3^4}{a_1^4}\right)$.
Therefore, the coefficient of  $x^{2^n-1}$ is zero if and only if
 $A(a_1,a_2,a_3)=0$, which
is equivalent to Equation~\eqref{eq:indicator-3}   by setting $d_1=a_2/a_1$ and
 $d_2=a_3/a_1$.

  In the case that the coefficient of $x^{2^n-1}$ vanishes, 
 the degree of the expression $x^{2^n-2}\prod_{s=1}^3(\tr_n(a_sx)+1)$ is   at most $n-1$, so $\deg(I_{|A})\le n-4$  and thanks to Corollary~\ref{cor:general-rest-deg-by}, $\deg(I_{|A})=n-4$.
\hfill$\Box$


\begin{remark} A natural question is to determine when the equation~\eqref{eq:indicator-3} in Theorem~\ref{Theo:inv-codim3} has solutions $(d_1,d_2)$ such that
$1,d_1,d_2$ are $\F_2$-linear independent (as $a_1,a_2,a_3$ in Theorem~\ref{Theo:inv-codim3} must be linearly independent). The field $\F_{2^n}$ must have
dimension at least 3 as a $\F_2$-vector space, i.e.\  $n\ge 3$. In fact, the equation~\eqref{eq:indicator-3} does have a solution in  $\Bbb{F}_{2^3}$
  by taking
 $d_1$ a primitive element of $\Bbb{F}_{2^3}$ such that $d_1^3=d_1+1$ and by taking  $d_2=d_1^2$. Therefore, whenever $n$ is a multiple of 3, the equation has solutions in $\F_{2^n}$, as $\F_{2^3}\subseteq\F_{2^n}$. The condition that $n$ is a multiple of 3 is sufficient, but not necessary for the existence of solutions, as we shall see in the next example.
\end{remark}

\begin{example}\label{ex:solutions-of-eq}
For $n=3,4,\ldots, 12$,   we computed the number of solutions
$(d_2,d_2)$ of the equation~\eqref{eq:indicator-3}, such that $(1,
d_1, d_2)$ are $\F_2$-linearly independent. This number, which we
will denote by $z(n)$, is equal to $24, 0, 0, 24, 168, 336, 528,840,
1848, 4224 $, respectively. The number of tuples $(a_1,a_2,a_3)$
such that  $d_1 = a_2/a_1$ and $d_2 = a_3/a_1$
would therefore be $z(n)(2^n-1)$. On the other hand, for each fixed linear space $V$ of codimension~3 there are  $(2^3-1)(2^3-2)(2^3-2^2)$ ways to express it as the solution set of the system of equations $\tr_n(a_ix)=0$, with $i\in \{1,2,3\}$.  Therefore, using Theorem~\ref{Theo:inv-codim3}, we see that the number of degree-drop spaces of codimension~3 on which the degree of the multiplicative inverse function $I$ drops by~3 can be computed as $\frac{(2^n-1)z(n)}{(2^3-1)(2^3-2)(2^3-2^2)}$, which gives $1,
0,
0,
9,
127,
510,
1606,
5115,
22517,
102960$ spaces, respectively; the ratio of these spaces out of all the
$\SqBinom{n}{3}_2$
linear spaces of codimension~3 is equal to $1$,
$0$,
$0$,
$0.00645$,
$0.01075$,
$0.00525$,
$0.00204$,
$0.00081$,
$0.00044$,
$0.00025$, respectively.  For $7\le n\le 12$ we noticed that this ratio is close to $\frac{1}{2^n}$. To conclude, it is quite rare that the degree of the multiplicative inverse function drops by~3 on a subspace of codimension~3, in most cases it drops only by~2.
\end{example}
\begin{example}
     For $n=8$, the multiplicative inverse function $I$ is affine equivalent to the S-Box of the AES cipher. For each linear space of codimensions $k\in \{1,2,3\}$ we computed the restriction of $I$ to that space, and the algebraic degree of that restriction.
      There are indeed no degree-drop space of codimension 1;
      on linear spaces of codimension~2, the degree drops by~1 on all but 85 spaces, on which it drops by~2.
      These 85 spaces are a proportion of $\frac{1}{2^7-1}\approx 0.00787$ of the 10795 $\F_2$-linear spaces of codimension~2 of $\F_{2^8}$ (and a proportion of $\frac{1}{2^8(2^7-1)}$ of all the affine spaces of codimension~2). The results of Theorem~\ref{thm:inv-codim-2} are confirmed.

    The degree drops by 2 on all the linear spaces of codimension~3 except for 510 spaces on which it drops by 3.
    This is a proportion of $\frac{510}{97155}\approx 0.00525$ of the total number of linear spaces of codimension~3, confirming the results of Theorem~\ref{Theo:inv-codim3}, combined with Example~\ref{ex:solutions-of-eq}.
\end{example}

\section{Counting the vectorial functions which have no degree-drop hyperplanes}\label{sec:cardofK1rnm}
In this section we will consider again $(n,m)$-functions for $n$ and $m$ positive and not necessarily equal; they will be represented in their multivariate ANF.


Denote by $e_1, \ldots, e_n$ the vectors of weight one, which form the canonical basis of $\F_2^n$.
For any $a \in \F_2^n\setminus\{\mathbf{0}\}$ we denote by $H_a$ the hyperplane defined by the equation $\sum_{i=1}^{n}a_ix_i=0$, which can also be written as $a \cdot x =0$ where ``$\cdot $'' denotes the usual scalar product.
From Lemma~\ref{lem:x1-factor} we have the obvious:
\begin{corollary}\label{cor:characterization}
  Let $F:\F_2^n \ra \F_2^m$ be a function of degree $r$. The hyperplane $H_{e_1}$ is a degree-drop hyperplane for $F$  if and only if all the monomials of degree $r$ of $F$ contain $x_1$.
\end{corollary}
It was shown in~\cite{SalMan17} that if $H_a$ and $H_b$  are degree-drop hyperplanes for a Boolean function $f$ then $H_{a+b}$ is also a degree-drop hyperplane for $f$. By applying these results to each coordinate of a vectorial function, we obtain the obvious generalization:
\begin{proposition}\label{vshyp}
  Let $F:\Bbb{F}_2^n \ra \Bbb{F}_2^m$ be a Boolean function
  and let $a,b \in \F_2^n\setminus\{\mathbf{0}\}$ with $a\neq b$. If $H_a$ and $H_b$  are degree-drop hyperplanes for $F$ then $H_{a+b}$ is also a degree-drop hyperplane for $F$.
\end{proposition}

Based on Proposition \ref{vshyp}, we can see that for each function
$F$ the set consisting of the zero vector and of all vectors $a\in
\F_2^n\setminus\{\mathbf{0}\}$ for which $H_a$ is a degree-drop hyperplane of
$F$ is a vector space.


Recall that we use $ \SqBinom{n}{k}_q$ to denote the  Gaussian
$q$-binomial coefficients, see~\eqref{eq:Gaussian-binom}.   We shall need  
the following result from \cite{SalFeru20}, which uses the result in
\cite{carlitz}:
\begin{lemma}\label{lem:Fourier-style}
Let $S,T: \mathbb{N} \rightarrow \mathbb{C}$ be functions.
Then
\begin{equation}\label{eq:rec-general}
  S(n) = \sum_{k=0}^{n}\SqBinom{n}{k}_q T(k) \mbox{  for all $n\ge 0$}
\end{equation}
if and only if
\begin{equation}\label{eq:explicit-general}
  T(n) = \sum_{k=0}^{n}(-1)^k q^{\frac{k(k-1)}{2}}\SqBinom{n}{k}_q S(n-k) \mbox{  for all $n\ge 0$}.
\end{equation}
\end{lemma}

We are now ready to compute the number of $(n,m)$-functions that do not have degree-drop hyperplanes.
\begin{theorem}\label{thm:V1-explicit-formula}
 The number of homogeneous $(n,m)$-functions of degree $r$
 which do not have any degree-drop hyperplane is equal to
 \begin{equation}\label{eq:K1-explicit-formula}
 \sum_{i=0}^{r}(-1)^i 2^{\frac{i(i-1)}{2}}\SqBinom{n}{i}_2 \left( 2^{m\binom{n-i}{r-i}}-1\right).
\end{equation}
  For any $j$ with $0\le j\le r$, the number of homogeneous $(n,m)$-functions of degree $r$ 
  which have exactly $2^j-1$ degree-drop linear hyperplanes  is equal to
   \begin{equation*}
 \SqBinom{n}{j}_2 \sum_{i=0}^{r-j}(-1)^i 2^{\frac{i(i-1)}{2}}\SqBinom{n-j}{i}_2 \left( 2^{m\binom{n-j-i}{r-j-i}}-1\right)
\end{equation*}
where  $\SqBinom{n}{k}_q$ denotes the Gaussian $q$-binomial coefficient.
\end{theorem}
{\em Proof.}
Let us denote by $K_{r,j,n,m}$ the set of homogeneous $(n,m)$-functions of degree $r$ which have exactly $2^j-1$ degree-drop linear hyperplanes. In particular, the functions in $K_{r,0,n,m}$ have no degree-drop hyperplane.

For any function $F$ consider the vector space $V_F$ consisting of the zero vector and of all vectors $a\in
\F_2^n\setminus\{0\}$ for which $H_a$ is a degree-drop hyperplane of
$F$.
Given a vector space $V$, let us denote by  $K_{r,V,n,m}$ the set of homogeneous  $(n,m)-$functions of degree $r$
for which $V_F$
equals $V$.

 For every affine automorphism $\varphi_M$ of $\F_2^n$, defined as $\varphi_M(x) = Mx$ with $M$ an invertible matrix, we know by Lemma \ref{lem:homomorphism} that $H_{b}$ is a degree-drop hyperplane for $F\circ \varphi_M$
 if and only if $b\in\varphi_{M^T}(V) \setminus\{0\}$.
We will  choose $M$ such that $\varphi_{M^T}(V) = E_{\dim(V)}$
where we denote by $E_i$ the space generated by the basis
$\{e_1,\ldots,e_i\}$. We have therefore
 \[
F\in K_{r,V,n,m}  \Leftrightarrow F\circ \varphi_M\in
K_{r,E_{\dim(V)},n,m}.
\]
Since $F_1\circ\varphi_{M^T} = F_2\circ\varphi_{M^T}$ if and only if
$F_1=F_2$, we have that
 \[
|K_{r,V,n,m} | = |K_{r,E_{\dim(V)},n,m}|.
\]
  In other words,
      the cardinality of $K_{r,V,n,m}$ only depends on the dimension of $V$ and not on the space $V$ itself.
%
%
%
%

Using
Lemma~\ref{lem:x1-factor},
we see that the functions  $F\in K_{r,E_i,n,m}$ are of the form $F(x_1,
\ldots, x_n) = x_1 x_2\ldots x_i G(x_{i+1}, \ldots, x_n)$, with $G$
a homogeneous $(n-i,m)$-function  of degree $r-i$ which does not have any degree-drop hyperplane. Therefore
\[
|K_{r,E_i,n,m} | = |K_{r-i,0,n-i,m}|.
\]
Using these equalities and  the Grassmannian $Gr_j(\F_2^n)$, which is  the set of all
$j$-dimensional linear subspaces of $\F_2^n$,  we have
\begin{eqnarray}
  |K_{r,j,n,m} | &=& |\bigcup_{V\in Gr_j(\F_2^n)}  K_{r,V,n,m}|  = \sum_{V\in Gr_j(\F_2^n)}|K_{r,V,n,m}|\nonumber\\
  &=& \SqBinom{n}{j}_2 |K_{r,E_j,n,m} | = \SqBinom{n}{j}_2 |K_{r-j,0,n-j,m}|.\label{eq:1}
\end{eqnarray}
We used the fact that for any two distinct spaces $V$
and $U$ the sets $K_{r,V,n,m}$ and $K_{r,U,n,m}$ are disjoint, and also the fact that $|Gr_j(\F_2^n)| = \SqBinom{n}{j}_2$.
The number of   $(n,m)$-homogeneous functions of degree $r$ is $2^{m\binom{n}{r}}-1$. Therefore
\begin{align}\label{eq:2}
2^{m\binom{n}{r}}-1  &= \sum_{j=0}^{n}|K_{r,j,n,m}| \nonumber\\& = \sum_{j=0}^{n}
\SqBinom{n}{j}_2 |K_{r-j,0,n-j,m} | \nonumber\\&= \sum_{j=n-r}^{n}
\SqBinom{n}{j}_2 |K_{r-n+j,0,j,m} |,
\end{align}
where in the last equation we replaced the index of summation $j$ by
$n-j$ and used the fact that $\SqBinom{n}{j}_2=\SqBinom{n}{n-j}_2$.
 We will now apply Lemma~\ref{lem:Fourier-style}. We will not use the variable $r$ but use a new variable $t= n-r$. Putting
 \[
 S(n) = 2^{m\binom{n}{r}}-1
 \]
 and
 \[
 T(n)  = |K_{n-t,0,n,m} |
 \]
 we see that~\eqref{eq:2} becomes~\eqref{eq:rec-general} from Lemma~\ref{lem:Fourier-style}. Therefore~\eqref{eq:explicit-general} in Lemma~\ref{lem:Fourier-style} must hold, which gives:
 \begin{eqnarray*}
  T(n)&=& |K_{n-t,0,n,m} | \\
  &=&
   \sum_{i=0}^{n}(-1)^i 2^{\frac{i(i-1)}{2}}\SqBinom{n}{i}_2 S(n-i) \\
   &=& \sum_{i=0}^{n}(-1)^i 2^{\frac{i(i-1)}{2}}\SqBinom{n}{i}_2  2^{m\binom{n-i}{n-i-t}}-1,
\end{eqnarray*}
which, after substituting $t=n-r$
gives the first formula in the
statement of the theorem. The second one follows from~\eqref{eq:1}.
\hfill $\Box$
\begin{remark}\label{rem:count-non-hom}
    If we want to count all functions $F$ of degree $r$ with no degree-drop hyperplanes (respectively $2^j-1$ linear degree-drop hyperplanes), regardless of whether $F$ is homogeneous or not, using Lemma~\ref{lem:homog}, we just need to take
    the number of homogeneous functions given by Theorem~\ref{thm:V1-explicit-formula} and multiply it by the total number of $(n,m)$-functions of degree at most $r-1$, i.e.\ $2^{m\sum_{d=0}^{r-1}\binom{n}{d}}$.
\end{remark}

We give now a  connection between functions having  no degree drop spaces and  functions which do not have ``fast
points" (defined below).  Recall that for a function $F$ and $a\in\F_2^n\setminus \{\mathbf{0}\}$, the
discrete derivative of  $F$ in the direction $a$ is defined as
$D_aF(x)=F(x+a)+F(x)$. It is known that $\deg(D_a F)\leq \deg(F)-1$ (see~\cite{Lai94} for the case of Boolean functions; for vectorial Boolean functions the result follows by noticing that if $F = (f_1,\ldots,f_m)$ then $D_a F = (D_a f_1,\ldots,D_a f_m) $). When the inequality is strict, $a$ is called a ``fast point'':

\begin{definition}\cite{Duan} Let $F:\mathbb{F}_{2}^n
\rightarrow\mathbb{F}_2^m$ be a non-constant vectorial Boolean function in n
variables. Any non-zero vector $a \in \mathbb{F}_{2}^n$ such that $\deg(D_aF)< \deg(F) - 1$  is called a {\em fast point} for $F$.
\end{definition}

The number of Boolean functions which have no fast points was given in~\cite{SalMan17}. We expand these results to vectorial Boolean functions (i.e. $m>1$).

For any monomial $t = x_{i_1}\cdots x_{i_r}$, the complement monomial $\prod_{i \in \{1,\dots n\}\setminus \{i_1, \ldots, i_r\}} x_i$ (which also equals $\frac{x_1\cdots x_n}{t}$) will be denoted  by $t^c$. For a homogeneous Boolean function $f= \sum_{i=1}^{\ell} t_{i}$ with $t_i$ monomials, we define $f^c= \sum_{i=1}^{\ell} (t_{i})^c$, with the convention that if $f$ is the identically zero function, then $f^c$ is the identically zero function. Finally,  for a homogeneous $(n,m)$-function $F=(f_1,...,f_m)$ we define $F^c=(f_1^c,...,f_m^c)$.
From \cite[Section 4]{Hou},  for any invertible matrix $M$  and any Boolean functions $f,g$, we have that $g = f \circ \varphi_{M}+h$ for some $h$ with $\deg(h)<\deg(f)$ if and only if $g^c = f^c \circ \varphi_{(M^T)^{-1}}+h'$ for some $h'$ with $\deg(h')<\deg(f^c)$. This can easily be generalized to vectorial functions.:
\begin{proposition}\label{prop:Hou-complement}
For any invertible matrix $M$ and any  $(n,m)$-functions $F,G$, we have that $G = F \circ \varphi_{M}+H$ for some $H$ with $\deg(H)<\deg(F)$ if and only if $G^c = F^c \circ \varphi_{(M^T)^{-1}}+H'$ for some $H'$ with $\deg(H')<\deg(F^c)$.
\end{proposition}
{\em Proof. }
Write $F=(f_1,\dots,f_m)$, $G=(g_1,\dots,g_m)$ and $H=(h_1,\dots,h_m)$. Then,  for all $i\in\{1,2\ldots m\}$, we have  $g_i= f_i \circ \varphi_{M}+h_i$ which is  equivalent to  $g_i^c = f_i^c \circ \varphi_{(M^T)^{-1}}+h_i'$ (for some $h_i'$ with $\deg(h_i')<\deg(f_i^c)$) that is, $G^c = F^c \circ \varphi_{(M^T)^{-1}}+H'$ where $H'=(h_1',\dots,h_m')$ with $\deg(H')<\deg(F^c)$.
\hfill$\Box$

We show that the existence of degree-drop hyperplanes for a homogeneous function $F$ is equivalent to the existence of fast points for the complement function $F^c$.
 We will actually prove a more general result, generalizing~\cite[Theorem~8]{CSA2}.

\begin{proposition}\label{thm:dd-space-fast-space}
  Let $F$ be a homogeneous vectorial function of algebraic degree $r$ in $n$ variables, with $r<n$. Let $1\le k\le n-r$, and let $a^{(1)}, \ldots, a^{(k)}$ be $k$ linearly independent elements of $\F_2^n$. The following statements are  equivalent:
\begin{itemize}
   \item The linear space defined by the $k$ equations
   $a^{(1)} \cdot x =0, \ldots, a^{(k)} \cdot x =0$
  is a  degree-drop subspace for $F$.
  \item
  Denoting $D^{(k)}_{a^{(1)}, \ldots, a^{(k)}}F = D_{a^{(1)}} (D_{a^{(2)}} (\ldots D_{a^{(k)}}F))$,  we have
  \[\deg(D^{(k)}_{a^{(1)}, \ldots, a^{(k)}} F^{c}) <\deg(F^c)-k,\]
  (i.e.\ the linear space generated by $a^{(1)}, \ldots, a^{(k)}$ is what is called a ``fast space'' for $F^c$ in~\cite[Definition~5]{SalMan17}, for $m=1$).
  \end{itemize}
\end{proposition}
 {\em Proof.} Let us set $F=(f_1,f_2,\ldots,f_m)$, where for all $i\in\{1,2,\ldots,m\}$ the Boolean function $f_i$ is either identically zero or homogeneous of degree $r$.
 For all $i\in\{1,2,\ldots,m\}$  for which $f_i$ is not the identically zero function, the linear space defined by the $k$ equations
  $a^{(1)} \cdot x =\ldots =a^{(k)} \cdot x =0$
  is a  degree-drop subspace for $f_i$ if and only if
  $\deg(D^{(k)}_{a^{(1)}, \ldots, a^{(k)}} f_i^{c}) <\deg(f_i^c)-k$, see \cite[Theorem~8]{CSA2}.
  The proof is completed by noting that $D^{(k)}_{a^{(1)}, \ldots, a^{(k)}} F^{c} = (D^{(k)}_{a^{(1)}, \ldots, a^{(k)}} f_1^{c}, \ldots, D^{(k)}_{a^{(1)}, \ldots, a^{(k)}} f_m^{c})$.
 \hfill $\Box$\\

For the particular case $k=1$, Proposition~\ref{thm:dd-space-fast-space} gives:
\begin{corollary}\label{cor:Vect-deg-drop-fast-point}
    Let $F$ be a homogeneous $(n,m)$-function of degree $r$
    with $1\le r\le n-1$ and let $a\in\F_2^n\setminus \{\mathbf{0}\}$. The linear hyperplane $H_a$ defined by the equation $a\cdot x = 0$ is a degree-drop hyperplane for $F$ if and only if $a$ is a fast point for $F^c$. Consequently, the number of linear degree-drop hyperplanes of $F$ is equal to the number of fast points of $F^c$.
\end{corollary}
Corollary~\ref{cor:Vect-deg-drop-fast-point} and Theorem~\ref{thm:V1-explicit-formula}, together with the fact that $\deg(F^c) = n - \deg(F)$ for any homogeneous $(n,m)$-function $F$, yield:
\begin{corollary}\label{cor:fast-points}
 Let $r,n,m$ be integers such that $1\le r \le n-1$.
 The number of homogeneous $(n,m)$-functions of degree $r$ 
 which do not have any fast point  is equal to
 \begin{equation}\label{eqf:K11-explicit-formula}
 \sum_{i=0}^{n-r}(-1)^i 2^{\frac{i(i-1)}{2}}\SqBinom{n}{i}_2 \left( 2^{m\binom{n-i}{n-r-i}}-1\right).
\end{equation}
  For any $j$ with $0\le j\le n-r$, the number of homogeneous $(n,m)$-functions of degree $r$ 
  which have exactly $2^j-1$ fast points  is equal to
   \begin{equation*}
 \SqBinom{n}{j}_2 \sum_{i=0}^{n-r-j}(-1)^i 2^{\frac{i(i-1)}{2}}\SqBinom{n-j}{i}_2 \left( 2^{m\binom{n-j-i}{n-r-j-i}}-1\right)
\end{equation*}
where  $\SqBinom{n}{k}_q$ denotes the Gaussian $q$-binomial coefficient.
\end{corollary}
\begin{example}
   In~\cite{KolPol24}, the authors study the sporadic Brinkmann-Leander-Edel-Pott function, which is a function $S:\F_{2^6}\rightarrow \F_{2^6}$ of degree 3, and is the only known APN function
  that is not equivalent to either a monomial or a quadratic function. They exploit the fact that $S$ has the maximum possible number of fast points for a function of degree $3$ in $6$ variables, namely $2^3-1$. This is a property that not many functions have; Corollary~\ref{cor:fast-points} quantifies how rare this property is. Namely, we can compute that only   85885 out of the $2^{120}-1$ homogeneous functions of degree 3 in 6 variables have $2^3-1$ fast points, that is a proportion of approximately $2^{-103}\approx 10^{-31}$ (the ratio being the same if we consider non-homogeneous functions, see Remark~\ref{rem:count-non-hom}).
\end{example}

\section{Conclusion}
In this paper,
we studied vectorial Boolean functions
which maintain their algebraic degree when restricted to any affine subspace of some codimension $k$; optimal functions $F$ are those that have this property for all $k$ up to and including $k=n-\deg(F)$, where $n$ is the number of variables. We showed that  affine injective functions, APN functions and power functions of the form $F(x)=x^{2^k-1}$ (including the multiplicative inverse function which is used in the AES S-box), or more generally $F(x)=x^{1+2^u+\ldots+2^{(r-1)u}}$ with $\gcd(u,n)=1$, are optimal from this point of view.
Finding other such optimal functions is an open problem.

For power functions, we gave necessary and sufficient conditions under which the degree stays unchanged when the function is restricted to subspaces of codimensions 1 and 2; we also gave sufficient conditions for arbitrary $k$; finding necessary and sufficient conditions for power functions for arbitrary $k$, as well as considering functions which are sums of two (or more) power functions, are open problems. For the multiplicative inverse function, we also showed that the degree does not decrease much on spaces of codimensions 2 and 3.
We mainly focussed in this paper on functions having the number of output bits $m$ equal to the number of input bits $n$,
and previous work~\cite{CSA} considered Boolean functions, i.e.\ $m=1$; the behaviour of vectorial Boolean functions for $1<m<n$ could also be explored in further work.

We also determined  a formula for the number of vectorial Boolean functions of a given degree $r$ in $n$ variables for which the degree does not change when the functions are restrited to hyperplanes. This also gives a formula for counting  vectorial Boolean functions of a given degree $r$ in $n$ variables which have no ``fast points'' with respect to differentiation.


\section*{Acknowledgment}
We thank Daniel Panario for his valuable insights on a result concerning irreducible polynomials (Proposition~\ref{prop:Irred-poly}).

We also extend our thanks to Marine Minier for recommending  relevant cryptanalysis articles that informed our work.

Finally, we acknowledge the financial support of the Engineering and Physical Sciences Research Council (EPSRC), UK, under grant EP/W03378X/1.
\bibliography{references}
\section{Appendix}
We give here an alternative proof of
Proposition~\ref{prop:Lin-indep}. We shall need the following
result, which might be known but we were unable to find a reference:
\begin{proposition}\label{prop:Irred-poly}Let $q$ be a prime power, let $g(x)$ be an irreducible polynomial in $\F_q[x]$ of
degree $u$ and let $n$ be co-prime with $u$. Then  $g(x)$ is
irreducible  over $\F_{q^n}$.
\end{proposition}
{\em Proof.} Since $g(x)$ is irreducible in $\F_q[x]$ and of degree
$u$, it splits in $\F_{q^u}$, that is, there is an element $\beta\in
\F_{q^u}$ such that $g(x)=\prod_{i=0}^{u-1}(x-\beta^{q^i})$.

Now, let $\gamma = \beta^{q^i}\in \F_{q^u}\subseteq \F_{q^{nu}}$ be
one of the roots of $g$. Let $g_1$ be the minimal polynomial of
$\gamma$ over $\F_{q^n}[x]$. By the definition of a minimal
polynomial, $g_1$ is irreducible over $\F_{q^n}$ and is a factor of
$g$; also, $\gamma^{q^n}, \gamma^{q^{2n}}, \ldots$ are also roots of
$g_1$. In particular,  $\gamma^{q^{n'n}}$ is also a root of $g_1$,
with $n'$ being the multiplicative inverse of $n$ modulo $u$, which
exists as $n$ and $u$ are coprime. We have $n'n = 1 +tu$ for some
integer $t$. Since $\gamma\in \F_{q^u}$, we have $\gamma^{q^u} =
\gamma$. Therefore, $\gamma^{q^{n'n}} = \gamma^{q^{1+tu}} =
\gamma^q$.
 So we proved that if $g_1$ has root $\gamma$ then it also has root $\gamma^q$. This also means $g_1$ has root $\gamma^{q^2}$ etc, so in fact all the roots of $g$ are roots of $g_1$. Since $g_1$ is a factor of $g$, we must have  $g = g_1$, so $g$ is irreducible over $\F_{q^n}[x]$.
\hfill$\Box$
\\
{\em Proof of Proposition~\ref{prop:Lin-indep}.}

If $a_1,a_2,\ldots,a_t$ are $\F_{q^u}$-linearly independent, then
they are also $\F_q$-linearly independent, as $\F_q\subseteq
\F_{q^u}$. Now let us assume  $a_1,a_2,\ldots,a_t$ are
$\F_q$-linearly independent. Since $a_1,a_2,\ldots,a_t\in \F_{q^n}$,
we can complete them with $a_{t+1},\ldots,a_n$ to a basis
$a_1,a_2,\ldots,a_n$ of $\F_{q^n}$ as a $\F_q$-vector space.

Consider  $\beta$ a primitive element of $\F_{q^u}$ and $g(x)\in
\F_q[x]$ its  corresponding primitive polynomial of degree $u$.
Clearly, $g(x)$ is irreducible in  $\F_q[x]$ and since
$\gcd(n,u)=1$, then, thanks to Proposition~\ref{prop:Irred-poly},
$g(x)$ is irreducible in $\F_{q^n}[x]$. We can therefore use $\beta$
to extend $\F_{q^n}$ to $\F_{q^{nu}}$, that is, for any element $y$
of $\F_{q^{nu}}$ there are uniquely determined elements $c_0,
\ldots, c_{u-1}\in \F_{q^n}$ such that $y = \sum_{j=0}^{u-1} c_j
\beta^j$. On the other hand, each $c_j$, being an element of
$\F_{q^n}$, can be written uniquely in base $a_1,a_2,\ldots,a_n$  as
$c_j = \sum_{i=1}^n c_{i,j}a_i$. This means that any element $y$ of
$\F_{q^{nu}}$ can be written as
\[y=\sum_{j=0}^{u-1}\sum_{i=1}^n c_{i,j}a_i\beta^j\]
 with $c_{i,j}\in \F_q$ uniquely determined. In other words,
the set  $B=\Big{\{}a_i\beta^j: \textrm{  }  (i,j)\in
\{1,\ldots,n\}\times\{0,\ldots,u-1\}\Big{\}}$ is a  basis of
$\F_{q^{nu}}$ over $\F_q$.

Now assume for a contradiction that there exist
$b_1,b_2,\ldots,b_t\in \F_{q^u}$ such that $\sum_{i=1}^tb_ia_i=0$
and $(b_1,b_2,\ldots,b_t)\neq 0$ (where $0$ here is the zero vector
of $\F_{q^u}^t$). Since $b_i\in \F_{q^u}$, we can write
$b_i=\sum_{j=0}^{u-1}b_{i,j}\beta^j$ with $b_{i,j}\in \F_q$. The
equation  $\sum_{i=1}^tb_ia_i=0$ becomes
$$\sum_{i=1}^t\sum_{j=0}^{u-1}b_{i,j}a_i\beta^j=0.$$
 Since $B$ is a basis of $\F_{q^{nu}}$ over $\F_q$, we have $b_{i,j}=0$ for all $(i,j)\in\{1,\ldots,t\}\times \{0,\ldots,u-1\}$. So, $b_i=0$ for all $i\in\{1,\ldots,t\}$, a contradiction.
\hfill$\Box$

\end{document}